\documentclass[12pt]{article}
\usepackage{amsfonts}
\usepackage{full page,
amssymb, amscd, amsmath, graphicx, %showkeys
}
\newtheorem{thm}{Theorem}[section]
\newtheorem{defn}[thm]{Definition}
\newtheorem{prop}[thm]{Proposition}
\newtheorem{cor}[thm]{Corollary}
\newtheorem{lemma}[thm]{Lemma}
\newtheorem{rema}[thm]{Remark}

\newcommand{\halmos}{\rule{1ex}{1.4ex}}

\newcommand{\nn}{\nonumber \\}

 \newcommand{\res}{\mbox{\rm Res}}

\renewcommand{\hom}{\mbox{\rm Hom}}

 \newcommand{\pf}{{\it Proof.}\hspace{2ex}}
 \newcommand{\epf}{\hspace*{\fill}\mbox{$\halmos$}}
 \newcommand{\epfv}{\hspace*{\fill}\mbox{$\halmos$}\vspace{1em}}

\newcommand{\mbar}{\Big\vert}
\newcommand{\lbar}{\bigg\vert}

\newcommand{\C}{\mathbb{C}}
\newcommand{\Z}{\mathbb{Z}}
\newcommand{\R}{\mathbb{R}}

\newcommand{\N}{\mathbb{N}}

\newcommand{\one}{\mathbf{1}}

\title{ {\bf Twist vertex operators for twisted modules} }
\date{}
\author{Yi-Zhi Huang}
\begin{document}

\bibliographystyle{alpha}
\maketitle
\begin{abstract}
We introduce and study twist vertex operators for a 
(lower-bounded generalized) twisted modules for a grading-restricted 
vertex (super)algebra.
We prove the duality property,  weak associativity, 
a Jacobi identity, a generalized commutator formula,
generalized weak commutativity, and convergence and commutativity for products of more than
two operators
involving twist vertex operators.
These properties of twist vertex operators
play an important role in the author's
recent general, direct and explicit construction of (lower-bounded generalized)
twisted modules.
\end{abstract}

\renewcommand{\theequation}{\thesection.\arabic{equation}}
\renewcommand{\thethm}{\thesection.\arabic{thm}}
\setcounter{equation}{0}
\setcounter{thm}{0}
\section{Introduction}

Frenkel, Lepowsky and Meurman 
introduced twisted modules associated to automorphisms of finite order
of a vertex operator algebra in their construction \cite{FLM1} \cite{FLM2}
\cite{FLM}
of the moonshine module vertex operator algebra $V^{\natural}$.  
In \cite{H-log-twisted-mod}, the author introduced 
(generalized) twisted module associated a general automorphism $g$ of 
 a vertex operator algebra $V$. Twisted modules for
vertex operator algebras  have been studied extensively in mathematics and 
physics. See for example, \cite{Le1}, \cite{FLM2}, 
\cite{Le2}, \cite{FLM}, \cite{D}, \cite{DL},
\cite{DonLM1}, \cite{DonLM2},  \cite{Li}, \cite{BDM}, \cite{DoyLM1},
\cite{DoyLM2}, \cite{BHL}, \cite{H-log-twisted-mod},
\cite{B}, \cite{Y}, and the references in these papers. 

Recently, the author has found a general and direct construction of
(lower-bounded generalized) twisted modules for a grading-restricted 
vertex (super)algebra $V$ \cite{H-const-twisted-mod}. In this construction, 
besides ``twisted fields'' or twisted vertex operators, one crucial ingredient in this
construction is what we call the ``twist fields" or ``twist vertex operators.'' 
Such an operators is a generalization of the vertex operator $Y_{WV}^{W}$ 
for a $V$-module $W$ given by the 
skew-symmetry in Section 5.6 in \cite{FHL} and 
is a special type of twisted intertwining 
operators obtained from twisted vertex operator for the (lower-bounded 
generalized) twisted module
using the skew-symmetry isomorphism $\Omega_{+}$ in \cite{H-twisted-int}.

The reason why twist vertex operators are important for the construction of 
twisted modules (and in particular, modules) is quite simple. To construct 
twisted modules (or even modules), we usually start with a set of generating twisted
fields corresponding to a set of generating fields of the algebra $V$. Usually it is not difficult to 
verify that these twisted fields satisfy weak commutativity or commutativity. But for twisted modules (or even modules), 
associativity is the main property to be verified and is not a consequence of the 
weak commutativity or commutativity. For a set 
of generating twisted fields,  we cannot 
formulate the associativity as one of our assumptions on these generating twisted fields.  But for a twisted module $W$, 
if we use the skew-symmetry isomorphism 
for intertwining operators to introduce an intertwining operator of type $\binom{W}{WV}$ called
twist vertex operator map, then the 
associativity for twisted vertex operators can be rewritten as 
the commutativity involving twisted vertex operators, twist vertex 
operators and vertex operators for the algebra $V$. 
Thus by introducing additional twist fields and assuming that the commutativity holds for the generating
twisted fields, twist fields and generating fields for $V$, we can construct twisted 
modules using the same method as in the first construction  of grading-restricted vertex algebras
in \cite{H-2-const}.

In this paper, we study twist vertex operators for a 
(lower-bounded generalized) twisted modules for a grading-restricted 
vertex (super)algebra $V$ associated to an automorphism $g$ of $V$.
We first prove that both the semisimple and unipotent parts of an automorphism $g$ of $V$ 
are also automorphisms of $V$ and the unipotent part of $g$ is in fact the exponential
of a derivation of $V$. Similar results also hold for twisted modules. 
We also prove or reformulate some results for twisted vertex operators, including 
weak commutativity, a Jacobi identity, a commutator formula, and 
convergence and commutativity for
products of more than two twisted vertex operators. 
We then prove the main results involving twist vertex operators, including duality,
weak associativity, a Jacobi identity, a generalized commutator formula, 
generalized weak commutativity, and convergence and commutativity for
products of more than two operators with one operator being a twist vertex operator.

In our formulations and proofs of the main results in this paper, we 
use both the complex variable approach and the formal variable approach. 
The formal variable approach can be used because we have explicit expressions 
of the correlation functions 
obtained by analytically extending the convergent products or iterates of twisted, untwisted 
and twist vertex operators. 
The convergence and commutativity are formulated and proved using 
complex variables approach while weak associativity, 
Jacobi identities, generalized commutator formulas and 
generalized weak commutativity are formulated 
using formal variable approach and proved using both approaches. 
The complex variable formulations often look simpler than the 
formal variable formulations because we can choose branches of multivalued correlation
functions and assume that the reader is familiar with multivalued 
analytic functions and their branches. On the other hand, although the formal variable formulations 
look more complicated, they reveal explicitly the information about suitable branches of multivalued 
correlation functions in terms of purely algebraic expressions. 

Our results are formulated and proved for lower-bounded generalized $g$-twisted $V$-modules. 
But these results hold for generalized $g$-twisted $V$-modules satisfying 
the lower-truncation properties in Remark \ref{lower-truncation}.

This paper is organized as follows:
In Section 2, we recall the notion of grading-restricted vertex (super)algebra and their 
variants. We also recall
some basic and useful properties of automorphisms of grading-restricted vertex algebras. 
In Section 3, we recall the notion of generalized twisted modules and their 
variants. We also recall and prove some properties of twisted vertex operators for 
lower-bounded generalized twisted modules. In Section 4, 
we introduce twist vertex operators for a lower-bounded 
twisted module. The main properties for these vertex operators are proved in this section.

\paragraph{Acknowledgments}
The author is grateful to Bin Gui for discussions on
the $L(-1)$-derivative property for twisted modules for vertex algebras.

\renewcommand{\theequation}{\thesection.\arabic{equation}}
\renewcommand{\thethm}{\thesection.\arabic{thm}}
\setcounter{equation}{0}
\setcounter{thm}{0}
\section{Grading-restricted vertex (super)algebras and automorphisms}

We recall the notions of grading-restricted vertex (super)algebra and (generalized) twisted module in this section. 
We also prove some properties of automorphisms of a grading-restricted vertex algebra and their actions on
a (generalized) twisted module. 

Before we recall the basic definitions, we first discuss the conventions used throughout 
the present paper. We shall use $z, z_{1}, z_{2}, \dots, \zeta, \zeta_{1}, \zeta_{2},
\xi, \xi_{1}, \xi_{2}, \dots $ to denote complex variables and 
$x, x_{1}, x_{2}, \dots, y, y_{1}, y_{2}, \dots$ to denote formal variables. Logarithms of formal
variables are also formal variables. 

For $z\in \C^{\times}$, we use $l_{p}(z)$ to denote $\log |z|+\arg z i+2p\pi i$, where 
$0\le \arg z<2\pi$. We shall always use expressions such as $(x_{1}-x_{2})^{n}$ for $n\in \C$ or 
$(\log (x_{1}-x_{2}))^{k}$ for $k\in \N$ as power series in $x_{2}$. 
For a semisimple operator $\mathcal{S}$, $x^{\mathcal{S}}$ or 
$(x_{1}-x_{2})^{\mathcal{S}}$ and so on
acting on the eigenspace for $\mathcal{S}$ with eigenvalue $\alpha\in \C$ are 
defined to be $x^{\alpha}$ or $(x_{1}-x_{2})^{\alpha}$ and so on. For a nilpotent operator 
$\mathcal{N}$, $x^{\mathcal{N}}$ or 
$(x_{1}-x_{2})^{\mathcal{N}}$ and so on are 
defined to be $e^{\mathcal{N}\log x}$ or $e^{\mathcal{N}\log (x_{1}-x_{2})}$ and so on. 

We first recall the definitions of grading-restricted vertex (super)algebra,
some variants and their automorphisms.

\begin{defn}\label{grvsa}
{\rm A {\it grading-restricted vertex superalgebra} is
a $\frac{\Z}{2}$-graded vector space
$V=\coprod_{n\in \frac{\Z}{2}}V_{(n)}$, 
equipped with a linear map
\begin{align*}
Y_{V}: V\otimes V&\to V[[x, x^{-1}]], \nn
u\otimes v&\mapsto Y_{V}(u, x)v,
\end{align*}
or equivalently, an analytic map
\begin{align*}
Y_{V}: \C^{\times}&\to \hom(V\otimes V, \overline{V}), \nn
z&\mapsto Y_{V}(\cdot, z)\cdot: u\otimes v\mapsto Y_{V}(u, z)v
\end{align*}
called the {\it vertex operator map}
and a {\it vacuum} $\mathbf{1}\in V_{(0)}$
satisfying the following axioms:

\begin{enumerate}

\item Axioms for the grading:
(a) {\it Grading-restriction condition}: When $n$ is sufficiently negative,
$V_{(n)}=0$ and $\dim V_{(n)}<\infty$ for $n\in \frac{\Z}{2}$.
(b) {\it $L(0)$-commutator formula}: Let $L_{V}(0): V\to V$
be defined by $L_{V}(0)v=nv$ for $v\in V_{(n)}$. Then
$$[L_{V}(0), Y_{V}(v, x)]=x\frac{d}{dx}Y_{V}(v, x)+Y_{V}(L_{V}(0)v, x)$$
for $v\in V$.

%\item {\it Lower-truncation condition for vertex operators}:
%For $u, v\in V$, $Y_{V}(u, x)v$ contain only finitely many negative
%power terms, that is, $Y_{V}(u, x)v\in V((x))$ (the space of formal
%Laurent series in $x$ with coefficients in $V$ and with finitely
%many negative power terms).

\item Axioms for the vacuum: (a) {\it Identity property}:
Let $1_{V}$ be the identity operator on $V$. Then
$Y_{V}(\mathbf{1}, x)=1_{V}$. (b)
{\it Creation property}: For $u\in V$,
$\lim_{x\to 0}Y_{V}(u, x)\mathbf{1}$ exists and is equal to $u$.

\item {\it $L(-1)$-derivative property} and {\it $L(-1)$-commutator formula}:
Let $L_{V}(-1): V\to V$ be the operator
given by
$$L_{V}(-1)v=\lim_{x\to 0}\frac{d}{dx}Y_{V}(v, x)\one$$
for $v\in V$. Then for $v\in V$,
$$\frac{d}{dx}Y_{V}(v, x)=Y_{V}(L_{V}(-1)v, x)=[L_{V}(-1), Y_{V}(v, x)].$$

\item {\it Duality}: For $u_{1}, u_{2}$ in either $\coprod_{n\in \Z}V_{(n)}$
or $\coprod_{n\in \Z+\frac{1}{2}}V_{(n)}$, $v\in V$ and
$v'\in V'$, the series
\begin{align}
&\langle v', Y_{V}(u_{1}, z_{1})Y_{V}(u_{2}, z_{2})v\rangle,\label{product12}\\
(-1)&^{|u_{1}| |u_{2}|}\langle v', Y_{V}(u_{2}, z_{2})Y_{V}(u_{1}, z_{1})v\rangle,
\label{product21}\\
\langle v&', Y_{V}(Y_{V}(u_{1}, z_{1}-z_{2})u_{2}, z_{2})v\rangle,\label{iterate}
\end{align}
where $|u_{1}|$ (or $|u_{2}|$) is $0$ if $u_{1}$ (or $u_{2}$) is in $\coprod_{n\in \Z}V_{(n)}$
and is $1$ if $u_{1}$ (or $u_{2}$) is in $\coprod_{n\in \Z+\frac{1}{2}}V_{(n)}$,
are absolutely convergent
in the regions $|z_{1}|>|z_{2}|>0$, $|z_{2}|>|z_{1}|>0$,
$|z_{2}|>|z_{1}-z_{2}|>0$, respectively, to a common rational function
in $z_{1}$ and $z_{2}$ with the only possible poles at $z_{1}, z_{2}=0$ and
$z_{1}=z_{2}$.

\end{enumerate}
In the case that $\coprod_{n\in \Z+\frac{1}{2}}V_{(n)}=0$, the grading-restricted vertex 
superalgebra just defined is called a {\it grading-restricted
vertex algebra}.}
\end{defn}

For simplicity, in the formulations, calculations, proofs and discussions of
the results in this paper, we always choose elements of $V$ to be either 
in $\coprod_{n\in \Z}V_{(n)}$
or $\coprod_{n\in \Z+\frac{1}{2}}V_{(n)}$ and we shall use $|v|$ to denote $0$ if $v\in \coprod_{n\in \Z}V_{(n)}$ 
and to denote $1$ if $v\in \coprod_{n\in \Z+\frac{1}{2}}V_{(n)}$. 

\begin{rema}\label{duality}
{\rm In Definition \ref{grvsa}, the duality property can be stated separately as 
three axioms, that is, the {\it rationality} (the convergence of (\ref{product12}), (\ref{product21})
and (\ref{iterate}) to rational functions in the regions $|z_{1}|>|z_{2}|>0$, $|z_{2}|>|z_{1}|>0$ 
and $|z_{2}|>|z_{1}-z_{2}|>0$, respectively), the {\it commutativity} (the statement that the rational 
functions to which (\ref{product12}) and (\ref{product21}) converge are equal) 
and the {\it associativity} (the statement that 
the (\ref{product12}) and (\ref{iterate}) are equal in the region $|z_{1}|>
|z_{2}|>|z_{1}-z_{2}|>0$). These axioms are not independent. In fact,
the associativity follows from the rationality and commutativity (see \cite{FHL}) and 
the commutativity also follows from the rationality and associativity (see \cite{H-jacobi-int-alg}).}
\end{rema}

\begin{defn}\label{qvsa}
{\rm A {\it quasi-vertex operator (super)algebra} or a {\it M\"{o}bius vertex algebra}
is a grading-restricted vertex (super)algebra
$(V, Y_{V},  \one)$ together with an operator $L_{V}(1)$ of weight $1$ on $V$ satisfying
\begin{eqnarray*}
[L_{V}(-1), L_{V}(1)]&=&-2L_{V}(0),\\
{[L_{V}(1), Y_{V} (v, x)]}&= &Y_{V} (L_{V}(1)v, x) + 2xY_{V} (L_{V}(0)v, x) 
+ x^{2}Y_{V} (L_{V}(-1)v, x)
\end{eqnarray*}
for $v\in V$. }
\end{defn}

%\begin{defn}\label{vosa}
%{\rm Let $(V, Y_{V},  \one)$ be a grading-restricted vertex (super)algebra.
%A {\it conformal element} of $V$ is an element $\omega\in V_{(2)}$
%satisfying the following axioms:
%
%\begin{enumerate}
%
%\item There exists $c\in \C$ such that
%$Y (\omega, x)\omega$ expanded as a $V$-valued Laurent series is equal
%to $L_{V}(-1)\omega x^{-1} + 2\omega x^{-2} + \frac{c}{2}\one x^{-4}$
%plus a $V$-valued power series in $x$.
%
%\item  $L_{V}(-1) = \res_{x}Y_{V}(\omega, x)$ and $L_{V}(0)=
%\res_{x}x Y_{V}(\omega, x)$ ($\res_{x}$ being the operation of taking the coefficient of $x^{-1}$
%of a Laurent series).
%
%\end{enumerate}
%A grading-restricted
%vertex (super)algebra equipped with a conformal element is called a {\it vertex operator
%(super)algebra} (or, more consistent with our terminology, a {\it grading-restricted 
%conformal vertex (super)algebra}). }
%\end{defn}

\begin{defn}
{\rm Let $V_{1}$ and $V_{2}$ be grading-restricted vertex superalgebras. A 
homomorphism from $V_{1}$ to $V_{2}$ is a weight-preserving linear map 
$g: V_{1}\to V_{2}$  such that $gY_{V_{1}}(u, x)v=Y_{V_{2}}(gu, x)gv$. 
An isomorphism from $V_{1}$ to $V_{2}$ is an invertible homomorphism 
from $V_{1}$ to $V_{2}$. When $V_{1}=V_{2}=V$, an isomorphism from $V$ to $V$ 
is called an automorphism of $V$. }
\end{defn}

In the rest of this paper, we fix a grading-restricted vertex superalgebra $V$ and 
an automorphism $g$ of $V$. We shall need the following linear algebra lemma:

\begin{lemma}\label{g-decomp}
There exist a semisimple operator $\mathcal{S}_{g}$ on $V$ with eigenvalues  of the form
$e^{2\pi i \alpha}$ for some $\alpha \in \C$ and an
operator $\mathcal{N}_{g}$ nilpotent on the homogeneous subspace $V_{(n)}$ 
of $V$ for $n\in \Z$
such that $\mathcal{S}_{g}$ and $\mathcal{N}_{g}$
commute with each other and 
$g=e^{2\pi i(\mathcal{S}_{g}+\mathcal{N}_{g})}$.
In particular, there exists an operator $\mathcal{L}_{g}=\mathcal{S}_{g}+\mathcal{N}_{g}$ 
on $V$ such that
$g=e^{2\pi i\mathcal{L}_{g}}$ and $\mathcal{S}_{g}$ and $\mathcal{N}_{g}$ are the semisimple
and nilpotent parts of $\mathcal{L}_{g}$. 
\end{lemma}
\pf
One proof can be found in Section 2 of \cite{HY}, where $\mathcal{S}_{g}$ 
and $\mathcal{N}_{g}$ are denoted by $\mathbf{A}_{V}$ and $\mathcal{N}$, respectively. 
\epfv

Let $P_{V}$ be set of $\alpha\in \C$ such that $0\le \Re(\alpha)<1$ and $e^{2\pi i \alpha}$ is an
eigenvalue of $g$. Then $V$ can be decomposed as a direct sum 
$$V=\coprod_{\alpha\in P_{V}}V^{[\alpha]},$$
where $V^{[\alpha]}$ for $\alpha\in P_{V}$ is the generalized eigenspace for $g$ 
with eigenvalue $e^{2\pi i\alpha}$. In particular, $V^{[\alpha]}$ is the eigenspace for $e^{2\pi\mathcal{S}_{g}}$
with eigenvalue $e^{2\pi i\alpha}$.
 
Apply Lemma 2.4 in \cite{HY} to $V$ viewed as a $1_{V}$-twisted 
$V$-module, we obtain:

\begin{prop}\label{semis-nilpt-isoms}
The operators $\mathcal{S}_{g}$ and $\mathcal{N}_{g}$ 
have the following properties:
\begin{enumerate}
\item The operators $e^{2\pi i\mathcal{S}_{g}}$ and $e^{2\pi i\mathcal{N}_{g}}$
are also automorphisms of $V$. 

\item The operator $\mathcal{N}_{g}$ is a derivation of $V$, that is, 
$$[\mathcal{N}_{g}, Y_{V}(u, x)]=Y_{V}(\mathcal{N}_{g}u, x).$$

\item For any $t\in \C$, $e^{t\mathcal{N}_{g}}$ is an automorphism of $V$. 
For any formal variable 
$y$ and $u, v\in V$, 
\begin{equation}\label{V-formal-N-g-conj}
e^{y\mathcal{N}_{g}}Y_{V}(u, x)v=Y_{V}(e^{y\mathcal{N}_{g}}u, x)e^{y\mathcal{N}_{g}}v.
\end{equation}
In particular, when $y=\log x_{0}$, 
\begin{equation}\label{V-formal-N-g-conj-spe}
x_{0}^{\mathcal{N}_{g}}Y_{V}(u, x)v=Y_{V}(x_{0}^{\mathcal{N}_{g}}u, x)
x_{0}^{\mathcal{N}_{g}}v.
\end{equation}
\end{enumerate}
\end{prop}

\begin{rema}
{\rm It is important to note that in general $\mathcal{S}_{g}$ might not be a derivation of $V$. 
The reason is that in general $\log (1-(e^{2\pi i \mathcal{S}_{g}}-1))$ expanded in nonnegative powers of 
$e^{2\pi i \mathcal{S}_{g}}-1$ might not be well defined. 
One can also see this from examples. For example, the automorphism group of 
the moonshine module vertex operator algebra $V^{\natural}$ is the Monster (see \cite{FLM}). 
So every automorphism $g$ of $V^{\natural}$ is of finite order and semisimple. In particular, 
$g=e^{2\pi i\mathcal{S}_{g}}$ for some operator $\mathcal{S}_{g}$. If $\mathcal{S}_{g}$
is a derivation of $V^{\natural}$, then for $t\in \C$, $e^{t\mathcal{S}_{g}}$ is also an automorphism of $V^{\natural}$ 
and thus there are infinitely many automorphisms of $V^{\natural}$. But the automorphism group 
of $V^{\natural}$ is the Monster, a finite simple group. Contradiction.}
\end{rema}

\renewcommand{\theequation}{\thesection.\arabic{equation}}
\renewcommand{\thethm}{\thesection.\arabic{thm}}
\setcounter{equation}{0}
\setcounter{thm}{0}
\section{Twisted modules and twisted vertex operators}

In this section,  we first recall various notions of (generalized) $g$-twisted $V$-module. 
Then we recall and prove some basic properties of twisted vertex operators 
for a generalized $g$-twisted $V$-module $W$.

\begin{defn}\label{defn-twisted-mod}
{\rm A {\it generalized  $g$-twisted 
$V$-module} is a ${\C}\times \Z_{2} \times \C/\Z$-graded
vector space 
$$W = \coprod_{n \in \C, s\in \Z_{2}, [\alpha]\in \C/\Z} W_{[n]}^{s; [\alpha]}
=\coprod_{n \in \C, s\in \Z_{2}, \alpha\in P_{W}} W_{[n]}^{s; [\alpha]}$$
 (graded by 
{\it weights}, {\it $\Z_{2}$-fermion number} and {\it $g$-weights}, 
where $P_{W}$ is the subset of 
the set $\{\alpha\in \C\;|\;\Re(\alpha)\in [0, 1)\}$ such that $W_{[n]}^{s; [\alpha]}\ne 0$
for $\alpha\in P_{W}$,
equipped with operators $L_{W}(0)$ and $L_{W}(-1)$ on $W$, a linear map
\begin{eqnarray*}
Y_{W}^g: V\otimes W &\to& W\{x\}[\log  x],\\
v \otimes w &\mapsto &Y_{W}^g(v, x)w=\sum_{n\in \C}\sum_{k\in \N}
(Y^{g}_{W})_{n, k}x^{n}(\log x)^{k}
\end{eqnarray*}
called {\it twisted vertex operator map} and an action of $g$ 
satisfying the following conditions:
\begin{enumerate}

\item The {\it equivariance property}: For $p \in \mathbb{Z}$, $z
\in \mathbb{C}^{\times}$,  $v \in V$ and $w \in W$, 
$$(Y^{g}_{W})^{p + 1}(gv,
z)w = (Y^{g}_{W})^{p}(v, z)w,$$
where for $p \in \mathbb{Z}$, $(Y^{g}_{W})^{ p}(v, z)$
is the $p$-th analytic branch of $Y_{W}^g(v, x)$.

\item The {\it identity property}: For $w \in W$, $Y_{W}^g({\bf 1}, x)w
= w$.

\item  The {\it duality property}: For
$u, v\in V$ (recall our convention that we always choose elements of $V$ 
to be either in $\coprod_{n\in \Z}V_{(n)}$ or $\coprod_{n\in \Z+\frac{1}{2}}V_{(n)}$
so that $|u|$ and $|v|$ are well defined), 
$w \in W$ and $w' \in W'$, there exists a
multivalued analytic function with preferred branch of the form
\begin{equation}\label{correl-fns}
f(z_1, z_2) = \sum_{i,
j, k, l = 0}^N a_{ijkl}z_1^{m_i}z_2^{n_j}({\rm log}z_1)^k({\rm
log}z_2)^l(z_1 - z_2)^{-t}
\end{equation}
for $N \in \mathbb{N}$, $m_1, \dots,
m_N$, $n_1, \dots, n_N \in \mathbb{C}$ and $t \in \mathbb{Z}_{+}$,
such that the series
\begin{align}
&\langle w', (Y^{g}_{W})^{ p}(u, z_1)(Y^{g}_{W})^{p}(v,
z_2)w\rangle,\label{prod12}\\
(-1)&^{|u||v|}\langle w', (Y^{g}_{W})^{ p}(v,
z_2)(Y^{g}_{W})^{p}(u, z_1)w\rangle ,\label{prod21}\\
&\langle w', (Y^{g}_{W})^{ p}(Y_{V}(u, z_1 - z_2)v,
z_2)w\rangle \label{iter}
\end{align}
are absolutely convergent in
the regions $|z_1| > |z_2| > 0$, $|z_2| > |z_1| > 0$, $|z_2| > |z_1
- z_2| > 0$, respectively, and their sums are equal to the branch
\begin{equation}\label{correl-fns-p-branch}
f^{p,p}(z_{1}, z_{2})= \sum_{i, j, k, l = 0}^N
a_{ijkl}e^{m_il_p(z_1)}e^{n_jl_p(z_2)}l_p(z_1)^kl_p(z_2)^l(z_1 -
z_2)^{-t}
\end{equation}
of $f(z_1, z_2)$ in the region $|z_1| > |z_2| > 0$, the region $|z_2| > |z_1| > 0$, 
the region given by $|z_2| > |z_1- z_2| > 0$ and $|\arg z_{1}-\arg z_{2}|<\frac{\pi}{2}$, respectively.

\item Properties about the gradings: (a) The  {\it $L(0)$-grading condition}:
For $w\in W_{[n]}$, $n\in \C$, there exists
$K\in \Z_{+}$ such that 
$(L_{W}(0)-n)^{K}w=0$. (b) The {\it $L(0)$-commutator formula}: 
$$[L_{W}(0), Y^{g}_{W}(v, z)]=z\frac{d}{dz}Y^{g}_{W}(v, z)+Y^{g}_{W}(L_{V}(0)v, z)$$
for $v\in V$. (c) The {\it $g$-grading condition}: For $\alpha\in P_{W}$,
$w \in W{[\alpha]}$, there exist $\Lambda\in \Z_{+}$ such that $(g-e^{2\pi \alpha i})^{\Lambda}w=0$. 
(d) The {\it $g$-compatibility condition}  and {\it $\Z_{2}$-fermion number compatibility condition}: 
For $u\in V$ and $w\in W$, 
$gY_{W}^{g}(u,x)w=Y_{W}^{g}(gu,x)gw$ and 
 $|Y_{W}^{g}(u,x)w|=|u|+|w|$.

\item The  {\it $L(-1)$-derivative property} and {\it $L(-1)$-commutator formula}: For $v \in V$,
\[
\frac{d}{dx}Y^g_{W}(v, x) =[L_{W}(-1), Y^{g}_{W}(v, x)].
\]

\end{enumerate}
A {\it lower-bounded  generalized $g$-twisted $V$-module} 
is a generalized $g$-twisted
$V$-module  $W$  such that
$W_{[n]} = 0$ when $\Re(n)<B$  for some $B\in \R$. A 
{\it grading-restricted generalized $g$-twisted $V$-module} or simply 
a {\it $g$-twisted $V$-module} is 
a lower-bounded generalized $g$-twisted $V$-module $W$ such that  for each $n \in
\mathbb{C}$, $\dim W_{[n]}<\infty$.  }
\end{defn}

For $w\in \coprod_{n\in \C, \alpha\in P_{W}}
W_{[n]}^{s; [\alpha]}$, we shall use $|w|$ to denote $s$. 
We shall also use the convention that in the formulations, calculations, proofs and 
discussions of the results in this paper, we always choose elements of $W$ to 
be in $\coprod_{n\in \C, \alpha\in P_{W}}
W_{[n]}^{s; [\alpha]}$ for either $s=0$ or $s=1$. Thus when we let $w\in W$, $|w|$ is always well defined.

We shall need the homogeneous subspaces $W_{[n]}=\coprod_{s=1, 2, \alpha\in P_{W}}
W_{[n]}^{s; [\alpha]}$ for $n\in \C$ and 
$W^{[\alpha]}=\coprod_{n\in \C, s=1, 2}
W_{[n]}^{s; [\alpha]}$ for $\alpha\in P_{W}$ of $W$.

%\begin{rema}
%{\rm Usually a $g$-twisted $V$-module is a grading-restricted generalized 
%$g$-twisted $V$-module $W$
%in the sense above such that the homogeneous subspaces 
%$W_{[n]}$ for $n\in \C$  $W$ are eigenspaces 
%for $L_{W}(0)$. In this paper, for simplicity,  we also call 
%grading-restricted generalized $g$-twisted $V$-module simply a 
%$g$-twisted $V$-module. 
%There should be no confusion caused by this simplification of terminology.}
%\end{rema}

\begin{rema}\label{lower-truncation}
{\rm For a lower-bounded 
generalized $g$-twisted $V$-module $W$, by (2.7) in \cite{HY} (or more explicitly,
the two displayed formula after (2.7) in \cite{HY}), for $u\in V^{[\alpha]}$,
\begin{equation}\label{lower-truncation-1}
Y^{g}_{W}(u,x)=\sum_{k=1}^{N}\sum_{n\in \alpha+\Z}
(Y^{g}_{W})_{n, k}(u)x^{-n-1}(\log x)^{k}.
\end{equation}
Moreover,
for $u\in V^{[\alpha]}$ and $w\in W$, 
$Y^{g}_{W}(u,
x)w$ has only finitely many terms containing $x^{-\alpha+n}$ for $n\in -\N$ 
and for $u\in V^{[\alpha]}$ and $w'\in W'$,  $\langle w', 
Y_{W}^{g}(u, x)\cdot\rangle$ has only 
finitely many terms containing
$x^{-\alpha+n}$ for $n\in \Z_{+}$.  
In the rest of the present paper, we shall study only lower-bounded 
generalized $g$-twisted $V$-modules though all the results still hold for 
generalized $g$-twisted $V$-modules such that the twisted vertex operators 
have these properties remarked above. }
\end{rema}

Since (\ref{prod12}), (\ref{prod21}) and (\ref{iter}) all converges absolutely
in the corresponding regions to the corresponding branches of the multivalued 
analytic function $f(z_1, z_2) $ with preferred branch, we shall denote
$f(z_1, z_2)$ by
$$F(\langle w', Y^{g}_{W}(u, z_1)Y^{g}_{W}(v,
z_2)w\rangle),$$
$$F((-1)^{|u||v|}\langle w', Y^{g}_{W}(v,
z_2)Y^{g}_{W}(u, z_1)w\rangle)$$
or 
$$F(\langle w', Y^{g}_{W}(Y_{V}(u, z_1 - z_2)v,
z_2)w\rangle).$$
We shall also denote the branch $f^{p,p}(z_{1}, z_{2})$ for $p\in \Z$ by
$$F^{p}(\langle w', Y^{g}_{W}(u, z_1)Y^{g}_{W}(v,
z_2)w\rangle),$$
$$F^{p}((-1)^{|u||v|}\langle w', Y^{g}_{W}(v,
z_2)Y^{g}_{W}(u, z_1)w\rangle)$$
or 
$$F^{p}(\langle w', Y^{g}_{W}(Y_{V}(u, z_1 - z_2)v,
z_2)w\rangle).$$
We shall use the similar notations to denote the multivalued analytic functions
or branches to which the products or iterates of more than two twisted vertex operators 
and vertex operators for the algebra converge. 

\begin{rema}\label{L(-1)-derivative}
{\rm The duality property in Definition \ref{defn-twisted-mod} can also be 
stated separately as three axioms, generalized rationality, commutativity and 
associativity. We do not state them
explicitly here since they are similar to the corresponding properties for $V$ 
(see Remark \ref{duality}).
But for $W$, associativity does not follow from 
the generalized rationality and commutativity while the commutativity does follow
from the generalized rationality and associativity. 
Also, it was observed by Bin Gui that the duality property for $W$ 
implies the $L(-1)$-derivative property
for $W$: For $v\in V$,
$$\frac{d}{dx}Y^g_{W}(v, x) =Y^{g}_{W}(L_{V}(-1)v, x).$$
In fact, for $w\in W$ and $w'\in W'$,
\begin{align}\label{L(-1)-derivative-1}
\langle w', (Y^{g}_{W})^{p}(u, z_1)w\rangle
&=F^{p}(\langle w', Y^{g}_{W}(u, z_1)Y^{g}_{W}(\one,
z_2)w\rangle)\nn
&=F^{p}(\langle w', Y^{g}_{W}(Y_{V}(u, z_1 - z_2)\one,
z_2)w\rangle)
\end{align}
for $u\in V$. Taking derivative with respect to $z_{1}$ on both sides of 
(\ref{L(-1)-derivative-1}),  using the $L(-1)$-derivative property for $V$ and 
using the duality property again, we obtain 
\begin{align}\label{L(-1)-derivative-2}
\left\langle w', \frac{d}{dz_{1}}(Y^{g}_{W})^{p}(u, z_1)w\right\rangle
&=F^{p}\left(\left\langle w', Y^{g}_{W}\left(\frac{\partial}{\partial z_{1}}
Y_{V}(u, z_1 - z_2)\one,
z_2\right)w\right\rangle\right)\nn
&=F^{p}(\langle w', Y^{g}_{W}(Y_{V}(L_{V}(-1)u, z_1 - z_2)\one,
z_2)w\rangle)\nn
&=\langle w', (Y^{g}_{W})^{p}(L_{V}(-1)u, z_1)w\rangle.
\end{align}
This is equivalent to the $L(-1)$-derivative property for $W$. 
}
\end{rema}

In the rest of this paper, we fix a lower-bounded $g$-twisted $V$-module $W$. 
It has been proved in Section 2 of \cite{HY} that the action 
of $g$ on $W$ has the following 
properties similar to those of the action of  $g$ on $V$:

\begin{lemma}
There exist a semisimple operator on $W$, denoted still by $\mathcal{S}_{g}$,  
such that for $\alpha\in P_{W}$, $W^{[\alpha]}$ is the eigenspace of $\mathcal{S}_{g}$
with the eigenvalue $e^{2\pi i\alpha}$ and an
operator nilpotent on any element of $W$ and preserving the 
gradings of $W$, denoted still by $\mathcal{N}_{g}$, 
such that $\mathcal{S}_{g}$ and $\mathcal{N}_{g}$ on $W$
commute with each other and 
$g=e^{2\pi i(\mathcal{S}_{g}+\mathcal{N}_{g})}$ on $W$.
In particular, there exists an operator $\mathcal{L}_{g}=\mathcal{S}_{g}+\mathcal{N}_{g}$ 
on $W$ such that
$g=e^{2\pi i\mathcal{L}_{g}}$ and $\mathcal{S}_{g}$ and $\mathcal{N}_{g}$ are the semisimple
and nilpotent parts of $\mathcal{L}_{g}$. 
\end{lemma}
\pf
Since $W$ is a direct sum of generalized eigenspaces of the action of $g$ on $W$,
we have  a simpler proof of the existence of $\mathcal{S}_{g}$, $\mathcal{N}_{g}$
and $\mathcal{L}_{g}$ than the one in \cite{HY} (in \cite{HY}, 
$\mathcal{S}_{g}$ and $\mathcal{N}_{g}$ are denoted by $\mathbf{A}_{W}$ 
and $\mathcal{N}$). Here we give this proof. 

For $w\in W^{[\alpha]}$, $\alpha\in P_{W}$, 
$(g-e^{2\pi \alpha})^{\Lambda}w=0$ for some $\Lambda\in \Z_{+}$. 
Define $\mathcal{S}_{g}$ on $W$ by $\mathcal{S}_{g}w=\alpha w$.
Then $e^{2\pi i \mathcal{S}_{g}}w=e^{2\pi i \alpha}w$ and 
\begin{align*}
(e^{-2\pi i \mathcal{S}_{g}}g-1_{W})^{\Lambda}w&=e^{-2\pi i \Lambda \mathcal{S}_{g}}(g
-e^{2\pi \alpha})^{\Lambda}w\nn
&=0.\end{align*}
So $e^{-2\pi i \mathcal{S}_{g}}g-1_{W}$ is nilpotent 
on $w$. Define
\begin{align}\label{defn-N-g}
\mathcal{N}_{g}w&=\frac{1}{2\pi i}\log (1_{W}+(e^{-2\pi i \mathcal{S}_{g}}g-1_{W}))w\nn
&=\sum_{j\in \Z_{+}}\frac{(-1)^{j+1}}{j}(e^{-2\pi i \mathcal{S}_{g}}g-1_{W})^{j}w
\end{align}
for $w\in W$.
Then $e^{-2\pi i \mathcal{S}_{g}}gw=e^{2\pi i\mathcal{N}_{g}}w$ and thus 
$gw=e^{2\pi i(\mathcal{S}_{g}+\mathcal{N}_{g})}w$ for $w\in W$. 
It is clear that $\mathcal{S}_{g}$
commutes with $g$ and thus with $\mathcal{N}_{g}$.
\epfv

\begin{prop}
For $u\in V$, $w\in W$, $t\in \C$ and formal variables $x, y$ and $x_{0}$, we have 
\begin{align}
e^{2\pi i\mathcal{S}_{g}}Y_{W}^{g}(u, x)w
&=Y_{W}^{g}(e^{2\pi i\mathcal{S}_{g}}u, x)e^{2\pi i\mathcal{S}_{g}}w,\label{S-g-comp}\\
e^{2\pi i\mathcal{N}_{g}}Y_{W}^{g}(u, x)w
&=Y_{W}^{g}(e^{2\pi i\mathcal{N}_{g}}u, x)e^{2\pi i\mathcal{N}_{g}}w, \label{N-g-comp}\\
\mathcal{N}_{g}Y_{W}^{g}(u, x)w&=Y_{W}^{g}(\mathcal{N}_{g}u, x)w
+Y_{W}^{g}(u, x)\mathcal{N}_{g}w,\label{N-g-derivayion}\\
e^{t\mathcal{N}_{g}}Y_{W}^{g}(u, x)w
&=Y_{W}^{g}(e^{t\mathcal{N}_{g}}u, x)e^{t\mathcal{N}_{g}}w,\\
e^{y\mathcal{N}_{g}}Y_{W}^{g}(u, x)w
&=Y_{W}^{g}(e^{y\mathcal{N}_{g}}u, x)e^{y\mathcal{N}_{g}}w, \label{formal-N-g-conj}\\
x_{0}^{\mathcal{N}_{g}}Y_{W}^{g}(u, x)w&=Y_{W}^{g}(x_{0}^{\mathcal{N}_{g}}u, x)x_{0}^{\mathcal{N}_{g}}w.\label{formal-N-g-conj-spe}
\end{align}
\end{prop}

By Lemma 2.3 in \cite{HY},  we have 
\begin{equation}\label{y-g-w-0}
Y^{g}_{W}(u, x)=(Y^{g}_{W})_{0}(x^{-\mathcal{N}_{g}}u, x),
\end{equation}
where 
$(Y^{g}_{W})_{0}(v, x)$  is the constant term of $Y^{g}_{W}(u, x)$ 
viewed as a power series in $\log x$. 
First we need another explicit form of the expansion of twisted vertex operators in the 
powers of the logarithm of the variable.

\begin{prop}
For $u\in V$, 
\begin{align}\label{y-g-w-0-conj}
Y^{g}_{W}(u, x)&=x^{-\mathcal{N}_{g}}(Y^{g}_{W})_{0}(u, x)
x^{\mathcal{N}_{g}}\nn
&=\sum_{k\in \N}\frac{(-1)^{k}}{k!}(\log x)^{k}
[\overbrace{\mathcal{N}_{g}, \cdots, [\mathcal{N}_{g}}^{k},
(Y^{g}_{W})_{0}(u, x)]\cdots],
\end{align}
\end{prop}
\pf
Using (\ref{formal-N-g-conj-spe}) with $x_{0}=x$ and 
(\ref{y-g-w-0}), we have 
\begin{align*}
Y^{g}_{W}(u, x)&=Y^{g}_{W}(x^{-\mathcal{N}_{g}}x^{\mathcal{N}_{g}}u, x)\nn
&=x^{-\mathcal{N}_{g}}Y^{g}_{W}(x^{\mathcal{N}_{g}}u, x)x^{\mathcal{N}_{g}}\nn
&=x^{-\mathcal{N}_{g}}(Y^{g}_{W})_{0}(u, x)
x^{\mathcal{N}_{g}}.
\end{align*}
This is the first equality in (\ref{y-g-w-0-conj}). Expanding 
$x^{-\mathcal{N}_{g}}(Y^{g}_{W})_{0}(u, x)
x^{\mathcal{N}_{g}}$ as a power series in $\log x$, we obtain the 
second equality in (\ref{y-g-w-0-conj}).
\epfv

We have the following
weak commutativity for twisted vertex operators:

\begin{prop}\label{locality}
For $u, v\in V$, let 
$M_{u, v}\in \Z_{+}$ such that $x^{M_{u, v}}Y_{V}(u, x)v\in V[[x]]$. Then
\begin{equation}\label{locality-f}
(x_{1}-x_{2})^{M_{u, v}}Y_{W}^{g}(u, x_{1})Y_{W}^{g}(v, x_{2})
=(x_{1}-x_{2})^{M_{u, v}}(-1)^{|u||v|}Y_{W}^{g}(v, x_{2})Y_{W}^{g}(u, x_{1}).
\end{equation}
\end{prop}
\pf
This weak commutativity can be derived easily from the Jacobi identity 
(\ref{jacobi-1}) below for 
twisted vertex operators. Here we directly derive it from 
the duality property so that the reader will be familiar with the complex variable 
approach which is necessary when we prove the convergence later. 

Since $x^{M_{u, v}}Y_{V}(u, x)v\in V[[x]]$, 
\begin{equation}\label{iter1}
(z_{1}-z_{2})^{M_{u, v}}\langle w', (Y^{g}_{W})^{ p}(Y_{V}(u, z_1 - z_2)v,
z_2)w\rangle
\end{equation} 
is a Laurent series in $z_{2}$ and $z_{1}-z_{2}$ with only nonnegative powers of $z_{1}-z_{2}$. 
By the duality property, (\ref{iter1}) is the expansion of 
\begin{equation}\label{correl-fns-p-branch-0}
\sum_{i, j, k, l = 0}^N
a_{ijkl}e^{m_il_p(z_1)}e^{n_jl_p(z_2)}l_p(z_1)^kl_p(z_2)^l (z_{1}-z_{2})^{M_{u, v}-t}.
\end{equation}
in the region given by $|z_{2}|>|z_{1}-z_{2}|>0$ and $|\arg z_{1}-\arg z_{2}|<\frac{1}{2}$
as a Laurent series of 
$z_{2}$ and $z_{1}-z_{2}$. Since (\ref{iter1}) has only nonnegative powers of $z_{1}-z_{2}$,
(\ref{correl-fns-p-branch-0}) cannot have a pole at $z_{1}-z_{2}=0$ and therefore 
can be rewritten as 
\begin{equation}\label{correl-fns-p-branch-1}
\sum_{i, j, k, l = 0}^{N'}
a'_{ijkl}e^{m_il_p(z_1)}e^{n_jl_p(z_2)}l_p(z_1)^kl_p(z_2)^l
\end{equation}
for some $N'\in \N$ and $a'_{ijkl}\in \C$.
Then from the duality property, for $w\in W$ and $w'\in W'$, 
\begin{equation}\label{prod12-1}
(z_{1}-z_{2})^{M_{u, v}}\langle w', (Y^{g}_{W})^{ p}(u, z_1)(Y^{g}_{W})^{p}(v,
z_2)w\rangle
\end{equation}
and 
\begin{equation}\label{prod21-1}
(z_{1}-z_{2})^{M_{u, v}}(-1)^{|u||v|}\langle w', (Y^{g}_{W})^{ p}(v,
z_2)(Y^{g}_{W})^{p}(u, z_1)w\rangle
\end{equation}
converges absolutely in the regions $|z_1| > |z_2| > 0$ and $|z_2| > |z_1| > 0$,
respectively, to (\ref{correl-fns-p-branch-1}).
But the expansions of (\ref{correl-fns-p-branch-1}) in the regions 
$|z_1| > |z_2| > 0$ and $|z_2| > |z_1| > 0$ are just
(\ref{correl-fns-p-branch-1})  itself. Therefore both
(\ref{prod12-1}) and (\ref{prod21-1}) are equal to the finite sum 
(\ref{correl-fns-p-branch-1}) in the region given by $z_{1}, z_{2}\ne 0$, $z_{1}\ne z_{2}$.
In particular, (\ref{prod12-1}) and (\ref{prod21-1}) are equal in this region. 

Since $w$ and $w'$ are arbitrary and $M_{u. v}$ are independent of $w$ and $w'$, we obtain 
$$(z_{1}-z_{2})^{M_{u. v}}(Y^{g}_{W})^{ p}(u, z_1)(Y^{g}_{W})^{p}(v,
z_2)=(z_{1}-z_{2})^{M_{u. v}}(-1)^{|u||v|} (Y^{g}_{W})^{ p}(v,
z_2)(Y^{g}_{W})^{p}(u, z_1)$$
as series of the form
$$\sum_{m, n\in \C, k, l\in \N}c_{m, n, k, l}e^{ml_{p}(z_{1})}
e^{nl_{p}(z_{2})}l_p(z_1)^kl_p(z_2)^l$$
for $c_{m, n, k, l}\in W$.
This is equivalent to (\ref{locality-f}).
\epfv

For the map $(Y_{W}^{g})_{0}$, there is a Jacobi identity obtained in \cite{B} 
and proved in \cite{HY} to be equivalent to the duality property for $Y^{g}_{W}$.
Here we reformulate it to obtain a Jacobi identity
 for the twisted vertex operator map $Y_{W}^{g}$.

\begin{thm}
For $u, v\in V$,
\begin{align}\label{jacobi-1}
x_0^{-1}&\delta\left(\frac{x_1 - x_2}{x_0}\right)
Y_{W}^{g}(u, x_1)
Y_{W}^{g}(v, x_2)- (-1)^{|u||v|} x_0^{-1}\delta\left(\frac{- x_2 + x_1}{x_0}\right)
Y_{W}^{g}(v, x_2)
Y_{W}^{g}(u, x_1)\nn
&= x_1^{-1}\delta\left(\frac{x_2+x_0}{x_1}\right)
Y_{W}^{g}\left(Y_{V}\left(\left(\frac{x_2+x_0}{x_1}\right)^{\mathcal{L}_{g}}
u, x_0\right)v, x_2\right).
\end{align}
\end{thm}
\pf
In the case that $u\in V^{[\alpha]}$, we have the Jacobi identity  
\begin{align}\label{jacobi}
x_0^{-1}\delta&\left(\frac{x_1 - x_2}{x_0}\right)(Y_{W}^{g})_{0}(u, x_1)
(Y_{W}^{g})_{0}(v, x_2)\nn
& \quad - (-1)^{|u||v|}x_0^{-1}\delta\left(\frac{- x_2 + x_1}{x_0}\right)
(Y_{W}^{g})_{0}(v, x_2)(Y_{W}^{g})_{0}(u, x_1)\nn
&= x_1^{-1}\delta\left(\frac{x_2+x_0}{x_1}\right)\left(\frac{x_2+x_0}{x_1}\right)^{\alpha}
(Y_{W}^{g})_{0}\left(Y_{V}\left(\left(1 + \frac{x_0}{x_2}\right)^{\mathcal{N}_{g}}u, x_0\right)v, x_2\right)
\end{align}
for $(Y_{W}^{g})_{0}$.
See \cite{HY} for the proof in the case that 
$V$ is a grading-restricted vertex algebra; in the case that $V$ is a grading-restricted 
vertex superalgebra, the proof is completely the same.
By (\ref{jacobi}) and (\ref{y-g-w-0}), we have 
\begin{align}\label{jacobi-.5}
x_0^{-1}&\delta\left(\frac{x_1 - x_2}{x_0}\right)
Y_{W}^{g}(u, x_1)
Y_{W}^{g}(v, x_2)w- (-1)^{|u||v|} x_0^{-1}\delta\left(\frac{- x_2 + x_1}{x_0}\right)
Y_{W}^{g}(v, x_2)
Y_{W}^{g}(u, x_1)w\nn
& = x_1^{-1}\delta\left(\frac{x_2+x_0}{x_1}\right)
\left(\frac{x_2+x_0}{x_1}\right)^{\alpha}
(Y_{W}^{g})_{0}\left(Y_{V}\left(\left(1 + \frac{x_0}{x_2}\right)^{\mathcal{N}_{g}}
x_{1}^{-\mathcal{N}_{g}}u, x_0\right)x_{2}^{-\mathcal{N}_{g}}v, x_2\right)w\nn
&= x_1^{-1}\delta\left(\frac{x_2+x_0}{x_1}\right)
\left(\frac{x_2+x_0}{x_1}\right)^{\alpha}
Y_{W}^{g}\left(Y_{V}\left(x_{2}^{\mathcal{N}_{g}}
\left(1 + \frac{x_0}{x_2}\right)^{\mathcal{N}_{g}}
x_{1}^{-\mathcal{N}_{g}}u, x_0\right)v, x_2\right)w\nn
&= x_1^{-1}\delta\left(\frac{x_2+x_0}{x_1}\right)
\left(\frac{x_2+x_0}{x_1}\right)^{\alpha}
Y_{W}^{g}\left(Y_{V}\left(\left(\frac{x_2+x_0}{x_1}\right)^{\mathcal{N}_{g}}
u, x_0\right)v, x_2\right)w.
\end{align}
The right-hand side of (\ref{jacobi-.5}) can be rewritten as 
the right-hand side of (\ref{jacobi-1}). So (\ref{jacobi-1}) holds. 
\epf

\begin{cor}
For $u, v\in V$, 
\begin{align}\label{commu-formu}
Y&_{W}^{g}(u, x_1)
Y_{W}^{g}(v, x_2)- (-1)^{|u||v|} 
Y_{W}^{g}(v, x_2)
Y_{W}^{g}(u, x_1)\nn
&= \res_{x_{0}} x_1^{-1}\delta\left(\frac{x_2+x_0}{x_1}\right)
Y_{W}^{g}\left(Y_{V}\left(\left(\frac{x_2+x_0}{x_1}\right)^{\mathcal{L}_{g}}
u, x_0\right)v, x_2\right).
\end{align}
\end{cor}

Next we prove that the product of $k$ twisted vertex operators is absolutely convergent
to a multivalued analytic function of a certain form.

\begin{thm}\label{k-prod-twted-vo}
For $v_{1}\in V^{[\alpha_{1}]}, \dots, 
v_{k}\in V^{[\alpha_{k}]}$, $w\in W$ and $w'\in W'$, 
the series
\begin{equation}\label{k-prod}
\langle w', (Y_{W}^{g})^{p}(v_{1}, z_{1})\cdots (Y_{W}^{g})^{p}(v_{k}, z_{k})w\rangle
\end{equation}
is absolutely convergent in the region $|z_{1}|>\cdots>|z_{k}|>0$. 
Moreover, there exists a 
multivalued analytic function with preferred branch of the form 
$$\sum_{n_{1}, \dots, n_{k}=0}^{N}f_{n_{1}\cdots n_{k+l}}(z_{1}, \dots,
z_{k})  
z_{1}^{-\alpha_{1}}\cdots z_{k}^{-\alpha_{k}}
(\log z_{1})^{n_{1}}\cdots (\log z_{k})^{n_{k}},$$
denoted by 
\begin{equation}\label{k-correl}
F(\langle w', Y_{W}^{g}(v_{1}, z_{1})\cdots Y_{W}^{g}(v_{k}, z_{k})w\rangle)
\end{equation}
using our notations introduced in the preceding section, 
where $N\in \N$ and $f_{i_{1}
\cdots i_{k}n_{1}\cdots n_{k}}(z_{1}, \dots,
z_{k})$ for $i_{1},
\dots, i_{k}$, $n_{1}, \cdots, n_{k}=0, \dots, N$ are rational functions
of $z_{1}, \dots, z_{k}$ with the only possible poles $z_{i}=0$ for
$i=1, \dots, k$, $z_{i}-z_{j}=0$
for $i, j=1, \dots, k$, $i\ne j$, such that the sum of 
the series (\ref{k-prod}) is equal to the branch
\begin{align}\label{k-correl-p}
F&^{p}(\langle w', Y_{W}^{g}(v_{1}, z_{1})\cdots Y_{W}^{g}(v_{k}, z_{k})w\rangle)\nn
&=\sum_{n_{1}, \dots, n_{k}=0}^{N}f_{n_{1}\cdots n_{k}}(z_{1}, \dots,
z_{k})e^{-\alpha_{1}l_{p}(z_{1})}\cdots e^{-\alpha_{k}l_{p}(z_{k})}
(l_{p}(z_{1}))^{n_{1}}\cdots (l_{p}(z_{k}))^{n_{k}},
\end{align}
of (\ref{k-correl})
in the region 
given by $|z_{1}|>\cdots>|z_{k}|>0$. 
In addition,  the orders of the pole $z_{i}=0$ of the rational functions 
$f_{n_{1}\cdots n_{k}}(z_{1}, \dots, z_{k})$ 
have a lower bound independent of $v_{q}$ for $q\ne i$ and $w'$;
the orders of the pole $z_{i}=z_{j}$ of the rational functions 
$f_{i_{1} \cdots i_{k}n_{1}\cdots n_{k}}(z_{1}, \dots, z_{k})$ 
have a lower bound independent of $v_{q}$ for $q\ne i, j$, 
$w$ and $w'$.
\end{thm}
\pf
Consider the formal series
\begin{equation}\label{convergence-1}
\prod_{i=l}^{k}x_{l}^{\alpha_{l}}\prod_{1\le i<j\le k}(x_{i}-x_{j})^{M_{v_{i}, v_{j}}}
\langle w', Y_{W}^{g}(v_{1}, x_{1})
\cdots Y_{W}^{g}(v_{k}, x_{k})w\rangle.
\end{equation}
For $1\le q\le k$, using (\ref{locality-f}), the series (\ref{convergence-1}) is equal to
\begin{align}\label{convergence-2}
&\prod_{i=l}^{k}x_{l}^{\alpha_{l}}
\prod_{1\le i<j\le k}(x_{i}-x_{j})^{M_{v_{i}, v_{j}}}(-1)^{|v_{q}||v_{q+1}|
+\cdots +|v_{q}||v_{k}|}\cdot\nn
&\quad \cdot 
\langle w', Y_{W}^{g}(v_{1}, x_{1})\cdots Y_{W}^{g}(v_{q-1}, x_{q-1})
Y_{W}^{g}(v_{q+1}, x_{q+1})\cdots Y_{W}^{g}(v_{k}, x_{k})
Y_{W}^{g}(v_{q}, x_{q})w\rangle.
\end{align}
Since by Remark \ref{lower-truncation}, 
$x_{q}^{\alpha_{q}}Y_{W}^{g}(v_{q}, x_{q})w$
is a Laurent series in $x^{q}$ with only finitely many negative power terms, 
(\ref{convergence-2}) is also a
Laurent series in $x_{q}$ with only finitely many negative power terms. 
So the same is true for (\ref{convergence-1}).  On the other hand, using (\ref{locality}) again,
(\ref{convergence-1}) is equal to
\begin{align}\label{convergence-3}
&\prod_{i=l}^{k}x_{l}^{\alpha_{l}}
\prod_{1\le i<j\le k}(x_{i}-x_{j})^{M_{v_{i}, v_{j}}}(-1)^{|v_{q}||v_{1}|
+\cdots +|v_{q}||v_{q-1}|}\cdot\nn
&\quad \cdot 
\langle w', 
Y_{W}^{g}(v_{q}, x_{q})Y_{W}^{g}(v_{1}, x_{1})\cdots Y_{W}^{g}(v_{q-1}, x_{q-1})
Y_{W}^{g}(v_{q+1}, x_{q+1})\cdots Y_{W}^{g}(v_{k}, x_{k})w\rangle.
\end{align}
Since by Remark \ref{lower-truncation}, $x_{q}^{\alpha_{q}}\langle w', 
Y_{W}^{g}(v_{q}, x_{q})\cdot\rangle$ is a Laurent series in $x_{q}$ with only finitely many positive power terms,
(\ref{convergence-3}) is also Laurent series in $x_{q}$ with only finitely many positive power terms.
So the same is true for (\ref{convergence-1}). Thus (\ref{convergence-1}) 
must be a Laurent polynomial in $x_{q}$ with  polynomials in $\log x_{q}$  as coefficients,
or equivalently, a  polynomial in $\log x_{q}$ with Laurent polynomials in $x_{q}$ as coefficients.
Since this is true for $q=1, \dots,  k$,  (\ref{convergence-1}) is  a polynomial in
$\log x_{1}, \dots, \log x_{k}$ with Laurent polynomials
in $x_{1}, \dots, x_{k}$  as coefficients.

On the other hand, by Remark \ref{lower-truncation}, we have
\begin{equation}\label{convergence-5}
\prod_{i=l}^{k}x_{l}^{-\alpha_{l}}\langle w', Y_{W}^{g}(v_{1}, x_{1})
\cdots Y_{W}^{g}(v_{k}, x_{k})w\rangle \in \C((x_{1}))\cdots ((x_{k})),
\end{equation}
where as usual, for a ring $R$ and a formal variable $x$, 
we use $R((x))$ to denote the ring of Laurent series in $x$ with coefficients 
in $R$ and finitely many negative power terms. 
Since $\C((x_{1}))\cdots ((x_{k}))$ is a ring and 
$\prod_{1\le i<j\le k}(x_{i}-x_{j})^{-M_{v_{i}, v_{j}}}$ is in fact in this ring,
$\prod_{1\le i<j\le k}(x_{i}-x_{j})^{M_{v_{i}, v_{j}}}$ is invertible in this 
ring with the inverse $\prod_{1\le i<j\le k}(x_{i}-x_{j})^{-M_{v_{i}, v_{j}}}$. 
Since (\ref{convergence-1}) is  a polynomial in
$\log x_{1}, \dots, \log x_{k}$ with Laurent polynomials
in $x_{1}, \dots, x_{k}$  as coefficients, it is also in this ring. 
Therefore  (\ref{convergence-5})  is equal to 
the product of $\prod_{1\le i<j\le k}(x_{i}-x_{j})^{-M_{v_{i}, v_{j}}}$ 
and a polynomial in
$\log x_{1}, \dots, \log x_{k}$ with Laurent polynomials
in $x_{1}, \dots, x_{k}$  as coefficients. 
So
$$\langle w', Y_{W}^{g}(v_{1}, x_{1})
\cdots Y_{W}^{g}(v_{k}, x_{k})w\rangle \in \C((x_{1}))\cdots ((x_{k}))$$
is equal to the product of $\prod_{i=l}^{k}x_{l}^{-\alpha_{l}}
\prod_{1\le i<j\le k}(x_{i}-x_{j})^{-M_{v_{i}, v_{j}}}$
and a polynomial in
$\log x_{1}, \dots, \log x_{k}$ with Laurent polynomials
in $x_{1}, \dots, x_{k}$  as coefficients. 
Thus (\ref{k-prod}) is absolutely convergent in the region $|z_{1}|>\cdots>|z_{k}|>0$
to an analytic function of the form (\ref{k-correl-p}).

The properties of the function (\ref{k-prod}) follow from Remark \ref{lower-truncation} 
and the duality property.
\epf

\begin{cor}
For $v_{1}, \dots,  v_{k}\in V$, $w\in W$, $w'\in W'$,
$p\in \Z$ and $\sigma\in S_{k}$, 
$$F^{p}(\langle w', Y_{W}^{g}(v_{1}, z_{1})\cdots Y_{W}^{g}(v_{k}, z_{k})w\rangle)
=\pm F^{p}(\langle w', Y_{W}^{g}(v_{\sigma(1)}, z_{\sigma(1)})\cdots 
Y_{W}^{g}(v_{\sigma(k)}, z_{\sigma(k)})w\rangle),$$
where the sign $\pm$ is uniquely determined by $|v_{1}, \dots, |v_{k}|$ and $\sigma$.
\end{cor}
\pf
This result follows immediately from 
Theorem \ref{k-prod-twted-vo} and the duality property.
\epfv

\renewcommand{\theequation}{\thesection.\arabic{equation}}
\renewcommand{\thethm}{\thesection.\arabic{thm}}
\setcounter{equation}{0}
\setcounter{thm}{0}
\section{Twist vertex operators}

We introduce twist vertex operators and  prove 
the main results on twist vertex operators in this section. 

Recall that $W$ is a lower-bounded generalized $g$-twisted $V$-module. 
We first introduce twist vertex operators. 
Let 
\begin{align*}
(Y^{g})_{WV}^{W}: W\otimes V&\to W\{x\}[\log  x],\nn
w \otimes v &\mapsto (Y^{g})_{WV}^{W}(w, x)v=\sum_{n\in \C}\sum_{k\in \N}
((Y^{g})_{WV}^{W})_{n, k}x^{n}(\log x)^{k}
\end{align*}
be defined by 
$$(Y^{g})_{WV}^{W}(w, x)v=(-1)^{|v||w|}e^{xL_{W}(-1)}Y_{W}^{g}(v, y)w
\mbar_{y^{n}=e^{\pi ni}x^{n},\; \log y=\log x+\pi i}$$
for $v\in V$ and $w\in W$. In particular, for $p\in \Z$, we have
\begin{align}\label{twist-vo}
((Y^{g})_{WV}^{W})^{p}(w, z)v&=(-1)^{|v||w|}e^{zL_{W}(-1)}Y_{W}^{g}(v, y)w
\mbar_{y^{n}=e^{n(l_{p}(z)+\pi i)},\; \log y=l_{p}(z)+\pi i}\nn
&=\left\{\begin{array}{ll}
(-1)^{|v||w|}e^{zL_{W}(-1)}(Y_{W}^{g})^{p}(v, -z)w&0\le \arg z< \pi,\\
(-1)^{|v||w|}e^{zL_{W}(-1)}(Y_{W}^{g})^{p+1}(v, -z)w&\pi \le \arg z< 2\pi.\end{array}\right.
\end{align}
We shall call $(Y^{g})_{WV}^{W}(w, x)$ a {\it twist vertex operator from $V$ to $W$} and 
$(Y^{g})_{WV}^{W}$ a {\it twist vertex operator map of type $\binom{W}{WV}$}.

Note that 
though the definitions and results in \cite{H-twisted-int} are given for a vertex operator 
algebras, they can be adapted to give definitions and results for grading-restricted 
vertex superalgebras by adding appropriate signs when two elements change their 
order. In particular, it is easy to see from the definition of generalized twisted module that
the twisted vertex operator map $Y^{g}_{W}$ is a twisted intertwining operator
of type $\binom{W}{VW}$.  By the definition of $(Y^{g})_{WV}^{W}$ above
and the definition of $\Omega_{+}$ in Section 5 of \cite{H-twisted-int}, 
we have
$$(Y^{g})_{WV}^{W}=\Omega_{+}(Y^{g}_{W}).$$
In fact there is another twisted intertwining operator $\Omega_{-}(Y^{g}_{W})$
of the same type. But we shall not use it in this paper. 
From the properties of $Y^{g}_{W}$ and Section 5 of \cite{H-twisted-int}, 
we obtain the following result:

\begin{thm}\label{twisted-io}
The linear map $(Y^{g})_{WV}^{W}$ is a twisted intertwining operator 
of type $\binom{W}{WV}$. In particular, we have the following 
properties of $(Y^{g})_{WV}^{W}$:
\begin{enumerate}

\item The lower truncation property: For $w\in W$ and $v\in V$,
$(Y^{g})_{WV}^{W}(w, x)v$ has only finitely many terms involving $x^{n}$ for $n\in \C$ 
with $\Re(n)<0$ and $(\log x)^{m}$ for $m\in \N$.

\item The duality property:  For $u, v\in V$, $w\in W$
and $w'\in W'$, the series 
\begin{align}
&\;\;\;\langle w', (Y_{W}
^{g})^{p_{1}}(u, z_{1})((Y^{g})_{WV}^{W})^{p_{2}}(w, z_{2})v\rangle,
\label{int-prod}\\
&\;(-1)^{|u||w|}\langle w', ((Y^{g})_{WV}^{W})^{p_{2}}(w, z_{2})
Y_{V}(u, z_{1})v\rangle,
\label{int-rev-prod}\\
 &\langle w', ((Y^{g})_{WV}^{W})^{p_{2}}(
(Y_{W}^{g})^{p_{1}}(u, z_{1}-z_{2})w, z_{2})v\rangle
\rangle\label{int-iter}
\end{align}
are absolutely convergent in the regions 
$|z_{1}|>|z_{2}|>0$,
$|z_{2}|>|z_{1}|>0$, $|z_{2}|>|z_{1}-z_{2}|>0$, respectively. Moreover, 
there exists a
multivalued analytic function with preferred branch
\begin{align}\label{f}
f(z_1, z_2; u, w, v, w')%\nn
&%&\quad 
= \sum_{j, k, m, n=0}^{N}a_{jkmn}z_1^{r}z_2^{s_{j}}(z_1 - z_2)^{t_{k}}
(\log z_2)^{m}(\log (z_1 - z_2))^{n}
\end{align}
for $N \in \mathbb{N}$, $r\in-\N$, $s_{j}, t_{k}, a_{ijklmn}\in \mathbb{C}$
such that for $p_{1}, p_{2}\in \Z$, their sums are equal 
to the branches
\begin{align}
f&^{p_{1}, p_{2}, p_{1}}(z_1, z_2; u, w, v, w')\nn
&=\sum_{j, k, m, n=0}^{N}a_{jkmn}z_{1}^{r}
e^{s_{j}l_{p_{2}}(z_2)}e^{t_{k}l_{p_{1}}(z_1 - z_2)}
(l_{p_{2}}(z_2))^{m}(l_{p_{1}}(z_1 - z_2))^{n},\label{f-prod}\\
f&^{p_{1}, p_{2}, p_{2}}(z_1, z_2; u, w, v, w')\nn
&=\sum_{j, k, m, n=0}^{N}a_{jkmn}z_{1}^{r}
e^{s_{j}l_{p_{2}}(z_2)}e^{t_{k}l_{p_{2}}(z_1 - z_2)}
(l_{p_{2}}(z_2))^{m}(l_{p_{2}}(z_1 - z_2))^{n},\label{f-rev-prod}\\
f&^{p_{2}, p_{2}, p_{1}}(z_1, z_2; u, w, v, w')\nn
&=\sum_{j, k, m, n=0}^{N}a_{jkmn}z_{1}^{r}
e^{s_{j}l_{p_{2}}(z_2)}e^{t_{k}l_{p_{1}}(z_1 - z_2)}
(l_{p_{2}}(z_2))^{m}(l_{p_{1}}(z_1 - z_2))^{n},\label{f-iter}
\end{align}
respectively, of $f(z_1, z_2; u, w_{1}, w_{2}, w_{3}')$ 
 in the region given by 
$|z_{1}|>|z_{2}|>0$ and $|\arg (z_{1}-z_{2})-\arg z_{1}|<\frac{\pi}{2}$,
the region given by $|z_{2}|>|z_{1}|>0$ and $-\frac{3\pi}{2}< \arg (z_{1}-z_{2})-\arg z_{2}<-\frac{\pi}{2}$,
the region given by $|z_{2}|>|z_{1}-z_{2}|>0$, respectively. In addition, when 
$u\in V^{[\alpha]}$ and $v\in V^{[\beta]}$, we can always take 
$a_{ijkl}=0$ for $j, k\ne 0$,  $s_{0}=-\beta$ and $t_{0}=-\alpha$.

\item The $L(-1)$-derivative property and $L(-1)$-commutator formula: 
\begin{align*}
\frac{d}{dx}(Y^{g})_{WV}^{W}(w, x)&=(Y^{g})_{WV}^{W}(L_{W}(-1)w, x)\nn
&=L_{W}(-1)(Y^{g})_{WV}^{W}(w, x)-(Y^{g})_{WV}^{W}(w, x)L_{V}(-1).
\end{align*}

\end{enumerate}
\end{thm}
\pf
Note that Theorem 5.1 in \cite{H-twisted-int} also holds for twisted intertwining operators 
for grading-restricted vertex superalgebras by adding appropriate signs.
By this theorem, the map $(Y^{g})_{WV}^{W}$ is a twisted intertwining operator 
of type $\binom{W}{WV}$. By the definition of twisted intertwining operator,
$(Y^{g})_{WV}^{W}$ have the lower truncation property,  the 
$L(-1)$-derivative property and the and $L(-1)$-commutator formula. 

By the duality property for twisted intertwining operators, 
(\ref{int-prod}), (\ref{int-rev-prod}) and (\ref{int-iter}) are absolutely convergent 
in the regions $|z_{1}|>|z_{2}|>0$,
$|z_{2}|>|z_{1}|>0$, $|z_{2}|>|z_{1}-z_{2}|>0$, respectively. Moreover 
there exists 
a multivalued analytic function 
\begin{align*}
f&(z_1, z_2; u, w, v, w')\nn
&\quad = \sum_{i, j, k, l, m, n=0}^{N}a_{ijklmn}z_1^{r_{i}}z_2^{s_{j}}(z_1 - z_2)^{t_{k}}(\log z_1)^l
(\log z_2)^{m}(\log (z_1 - z_2))^{n}
\end{align*}
for $N \in \mathbb{N}$, $r_{i}, s_{j}, t_{k}, a_{ijklmn}\in \mathbb{C}$,
such that
(\ref{int-prod}), (\ref{int-rev-prod}) and (\ref{int-iter}) are equal 
to the branches
\begin{align*}
f&^{p_{1}, p_{2}, p_{1}}(z_1, z_2; u, w, v, w')\nn
&=\sum_{i, j, k, l, m, n=0}^{N}a_{ijklmn}e^{r_{i}l_{p_{1}}(z_{1})}
e^{s_{j}l_{p_{2}}(z_2)}e^{t_{k}l_{p_{1}}(z_1 - z_2)}(l_{p_{1}}(z_1))^l
(l_{p_{2}}(z_2))^{m}(l_{p_{1}}(z_1 - z_2))^{n},\nn
f&^{p_{1}, p_{2}, p_{2}}(z_1, z_2; u, w, v, w')\nn
&=\sum_{i, j, k, l, m, n=0}^{N}a_{ijklmn}e^{r_{i}l_{p_{1}}(z_{1})}
e^{s_{j}l_{p_{2}}(z_2)}e^{t_{k}l_{p_{2}}(z_1 - z_2)}(l_{p_{1}}(z_1))^l
(l_{p_{2}}(z_2))^{m}(l_{p_{2}}(z_1 - z_2))^{n},\nn
f&^{p_{2}, p_{2}, p_{1}}(z_1, z_2; u, w, v, w')\nn
&=\sum_{i, j, k, l, m, n=0}^{N}a_{ijklmn}e^{r_{i}l_{p_{2}}(z_{1})}
e^{s_{j}l_{p_{2}}(z_2)}e^{t_{k}l_{p_{1}}(z_1 - z_2)}(l_{p_{2}}(z_1))^l
(l_{p_{2}}(z_2))^{m}(l_{p_{1}}(z_1 - z_2))^{n},
\end{align*}
respectively, of $f(z_1, z_2; u, w, v, w')$ in the region given by 
$|z_{1}|>|z_{2}|>0$ and $|\arg (z_{1}-z_{2})-\arg z_{1}|<\frac{\pi}{2}$,
the region given by $|z_{2}|>|z_{1}|>0$ and $-\frac{3\pi}{2}<
 \arg (z_{1}-z_{2})-\arg z_{2}<-\frac{\pi}{2}$,
the region given by $|z_{2}|>|z_{1}-z_{2}|>0$ and $|\arg  z_{1}- \arg z_{2}|<\frac{\pi}{2}$, respectively. Since (\ref{int-rev-prod}) contains only integral powers of $z_{1}$,
$f^{p_{1}, p_{2}, p_{2}}(z_1, z_2; u, w, v, w')$ must be independent of $p_{1}$ 
and of the form of the right-hand
side of (\ref{f-rev-prod}). Thus $f(z_1, z_2; u, w, v, w')$ must be of the form 
of  (\ref{f}). In particular, 
(\ref{int-prod}) and (\ref{int-iter}) must be of the forms of
(\ref{f-prod}) and (\ref{f-iter}), respectively. 

When $u\in V^{[\alpha]}$, by (\ref{lower-truncation-1}), 
$$(Y_{W}^{g})^{p_{1}}(u, z_{1}-z_{2})w=\sum_{k=0}^{N}\sum_{n\in \alpha+\Z}
(Y_{W}^{g})_{n, k}(u)we^{(-n-1)l_{p_{1}}(z_{1}-z_{2})}(l_{p_{1}}(z_{1}-z_{2}))^{k}.$$
When $v\in V^{[\beta]}$,
by the definition of $(Y^{g})_{WV}^{W}$ and (\ref{lower-truncation-1}),
we have 
\begin{align*}
&((Y^{g})_{WV}^{W})^{p_{2}}(w, z_{2})v\nn
&\quad =(-1)^{|v||w|}e^{z_{2}L_{W}(-1)}Y_{W}^{g}(v, y)w
\mbar_{y^{n}=e^{n(l_{p_{2}}(z_{2})+\pi i)},\; \log y=l_{p_{2}}(z_{2})+\pi i}\nn
&\quad =(-1)^{|v||w|}e^{z_{2}L_{W}(-1)}\sum_{l=0}^{M}\sum_{m\in \beta+\Z}
(Y_{W}^{g})_{m, l}(v)w e^{(-m-1)(l_{p_{2}}(z_{2})+\pi i)}(l_{p_{2}}(z_{2})+\pi i)^{l}
\end{align*}
Since (\ref{int-iter}) and (\ref{int-prod}) are convergent absolutely to 
$f^{p_{2}, p_{2}, p_{1}}(z_1, z_2; u, w, v, w')$ in the region given by 
$|z_{2}|>|z_{1}-z_{2}|>0$ and $|z_{1}|>|z_{2}|>0$, respectively, 
we can always choose the coefficients
$a_{ijkl}$ for $i, j, k, l=0, \dots, N$ such that 
$a_{ijkl}=0$ for $j, k\ne 0$, $s_{0}=-\beta$ and $t_{0}=-\alpha$.
\epfv

\begin{rema}
{\rm Using the definition of $(Y^{g})_{WV}^{W}$ and Remark \ref{lower-truncation}
we have the following stronger property than the lower-truncation property:
For $w\in W$ and $v\in V^{[\alpha]}$, there exist $N\in \N$ and $M\in \Z$ such that
$$(Y^{g})_{WV}^{W}(w, x)v=\sum_{k=1}^{N}\sum_{n\in \alpha+M-\N}
((Y^{g})_{WV}^{W})_{n, k}(w)v x^{-n-1}(\log x)^{k}.$$
For $w'\in W'$ and $w\in W$, 
$\langle w', 
(Y^{g})_{WV}^{W}(w, x)\cdot\rangle$ has only finitely many 
terms involving $x^{n}$ for $n\in \C$ 
with $\Re(n)>0$ and $(\log x)^{m}$ for $m\in \N$.}
\end{rema}

We shall use our notations introduced in Section 2 to 
denote the multivalued analytic function $f(z_1, z_2; u, w, v, w')$ with 
preferred branch in Theorem 
\ref{twisted-io}  by 
$$F(\langle w', Y_{W}
^{g}(u, z_{1})(Y^{g})_{WV}^{W}(w, z_{2})v\rangle)$$
or
$$F((-1)^{|u||w|}\langle w', (Y^{g})_{WV}^{W}(w, z_{2})
Y_{V}(u, z_{1})v\rangle)$$
or
$$F(\langle w', (Y^{g})_{WV}^{W}(
Y_{W}^{g}(u, z_{1}-z_{2})w, z_{2})v\rangle).$$
We shall also denote the branches (\ref{f-prod}), (\ref{f-rev-prod}),
(\ref{f-iter})  of $f(z_1, z_2; u, w, v, w')$ when $p_{1}=p_{2}=p$ by
\begin{align*}
&\quad\quad  F^{p}(\langle w', Y_{W}
^{g}(u, z_{1})(Y^{g})_{WV}^{W}(w, z_{2})v\rangle),\\
&F^{p}((-1)^{|u||w|}\langle w', (Y^{g})_{WV}^{W}(w, z_{2})
Y_{V}(u, z_{1})v\rangle),\\
&\quad F^{p}(\langle w', (Y^{g})_{WV}^{W}(
Y_{W}^{g}(u, z_{1}-z_{2})w, z_{2})v\rangle),
\end{align*}
respectively. Then in particular, we have the following commutativity 
for twisted and twist vertex operators:

\begin{cor}\label{twisted-twist-comm}
For $u, v\in V$, $w\in W$
and $w'\in W'$, we have 
\begin{equation}\label{twted-twt-coomm-f}
F^{p}(\langle w', Y_{W}^{g}(u, z_{1})(Y^{g})_{WV}^{W}
(w, z_{2})v\rangle)
=(-1)^{|u||w|} F^{p}(\langle w',(Y^{g})_{WV}^{W}
(w, z_{2}) Y_{V}(u, z_{1})v\rangle).
\end{equation}
\end{cor}

We also have the following formal weak associativity property involving twist vertex operators:

\begin{prop}\label{wk-assoc}
For $u, v\in V$, 
\begin{equation}\label{wk-assoc-1}
(x_{0}+x_{2})^{M_{u, v}}Y_{W}^{g}(u, x_{0}+x_{2})(Y^{g})_{WV}^{W}(w, x_{2})v
=(x_{0}+x_{2})^{M_{u, v}}(Y^{g})_{WV}^{W}(Y_{W}^{g}(u, x_{0})w, x_{2})v,
\end{equation}
where $M_{u, v}$ is the positive integer in Proposition \ref{locality}
\end{prop}
\pf
We take $u, v\in V$ and 
$w\in W$. Then by the definition of $(Y^{g})_{WV}^{W}$,
the $L(-1)$-commutator and $L(-1)$-derivative properties,  Proposition \ref{locality}
and the $\Z_{2}$-fermion number compatibility condition,
\begin{align*}
&(x_{0}+x_{2})^{M_{u, v}}Y^{g}_{W}(u, x_0+x_{2})
(Y^{g})_{WV}^{W}(w, x_{2})v\nn
&\quad =(x_{0}+x_{2})^{M_{u, v}} Y^{g}_{W}(u, x_0+x_{2})
(-1)^{|v||w|}Y_{W}^{g}(v, y)w
\mbar_{y^{n}=e^{\pi ni}x_{2}^{n},\; \log y=\log x_{2}+\pi i}\nn
&\quad =(-1)^{|v||w|}e^{x_{2}L_{W}(-1)}(x_{0}-y)^{M_{u, v}}
Y^{g}_{W}(u, x_0)
Y_{W}^{g}(v, y)w\mbar_{y^{n}=e^{\pi ni}x_{2}^{n},\; \log y=\log x_{2}+\pi i}
\nn
&\quad =(-1)^{|v||w|}(-1)^{|v||u|}e^{x_{2}L_{W}(-1)}(x_{0}-y)^{M_{u, v}}
Y_{W}^{g}(v, y)
Y^{g}_{W}(u, x_0)w\mbar_{y^{n}=e^{\pi ni}x_{2}^{n},\; \log y=\log x_{2}+\pi i}\nn
&\quad =(x_{0}-y)^{M_{u, v}}(-1)^{|v||Y^{g}_{W}(u, x_0)w|}e^{x_{2}L_{W}(-1)}
Y_{W}^{g}(v, y)
Y^{g}_{W}(u, x_0)w\mbar_{y^{n}=e^{\pi ni}x_{2}^{n},\; \log y=\log x_{2}+\pi i}\nn
&\quad =(x_{0}+x_{2})^{M_{u, v}}
(Y^{g})_{WV}^{W}(Y^{g}_{W}(u, x_0)w, x_{2})v.
\end{align*}
\epfv

Recall  (\ref{y-g-w-0}) (see Lemma 2.3 in \cite{HY}). 
We have a similar
result for $(Y^{g})_{WV}^{W}$. 
For $w\in W$, let 
$$((Y^{g})_{WV}^{W})_{0}(w, x)
=(Y^{g})_{WV}^{W}(w, x)x^{\mathcal{N}_{g}}.$$
Then we have: 

\begin{lemma}
For $w\in W$, $((Y^{g})_{WV}^{W})_{0}(w, x)v\in W\{x\}$
and 
\begin{equation}\label{twist-vo-decomp-2}
(Y^{g})_{WV}^{W}(w, x)
=((Y^{g})_{WV}^{W})_{0}(w, x)x^{-\mathcal{N}_{g}}.
\end{equation}
In particular, $((Y^{g})_{WV}^{W})_{0}(w, x)$ is the constant term of
$(Y^{g})_{WV}^{W}(w, x)$ viewed as a power series in $\log x$.
\end{lemma}
\pf 
By definition, for $v\in V$ and $w\in W$,
\begin{align}\label{twist-vo-decomp}
(Y^{g})_{WV}^{W}(w, x)v&=(-1)^{|v||w|}e^{xL_{W}(-1)}Y_{W}^{g}(v, y)w
\mbar_{y^{n}=e^{\pi ni}x^{n},\; \log y=\log x+\pi i}\nn
&=(-1)^{|v||w|}e^{xL_{W}(-1)}(Y_{W}^{g})_{0}(y^{-\mathcal{N}_{g}}v, y)w
\mbar_{y^{n}=e^{\pi ni}x^{n},\; \log y=\log x+\pi i}\nn
&=(-1)^{|v||w|}e^{xL_{W}(-1)}(Y_{W}^{g})_{0}(e^{-\log y\mathcal{N}_{g}}v, y)w
\mbar_{y^{n}=e^{\pi ni}x^{n},\; \log y=\log x+\pi i}\nn
&=(-1)^{|v||w|}e^{xL_{W}(-1)}(Y_{W}^{g})_{0}(x^{-\mathcal{N}_{g}}
e^{-\pi i\mathcal{N}_{g}}v, y)w
\mbar_{y^{n}=e^{\pi ni}x^{n}}.
\end{align}
Replacing $v$ in (\ref{twist-vo-decomp}) by $x^{\mathcal{N}_{g}}$,
we obtain
\begin{equation}\label{twist-vo-decomp-1}
(Y^{g})_{WV}^{W}(w, x)x^{\mathcal{N}_{g}}v
=(-1)^{|v||w|}e^{xL_{W}(-1)}(Y_{W}^{g})_{0}(e^{-\pi i\mathcal{N}_{g}}v, y)w
\mbar_{y^{n}=e^{\pi ni}x^{n}}.
\end{equation} 
From (\ref{twist-vo-decomp-1}) and the fact that 
$(Y_{W}^{g})_{0}(v, x)w\in W\{x\}$, we see that 
$((Y^{g})_{WV}^{W})_{0}(w, x)v\in W\{x\}$ and  (\ref{twist-vo-decomp-2}) holds.
\epfv

We now prove a Jacobi identity involving $(Y^{g})_{WV}^{W}$. 
In the results below, for simplicity, 
we shall always use the convention that for any operator 
or number $A$, 
\begin{equation}\label{convention}
(-x_{2}+x_1)^{A}=(x_{2}-x_1)^{A}e^{\pi iA}.
\end{equation}

\begin{thm}
For $u, v\in V$ and $w\in W$, we have
\begin{align}\label{jacobi-twist}
x_0&^{-1}\delta\left(\frac{x_1- x_2}{x_0}\right)
Y_{W}^{g}\left(\left(\frac{x_{1}-x_2}{x_0}\right)
^{\mathcal{L}_{g}}u, x_{1}\right)
(Y^{g})_{WV}^{W}(w, x_{2})v\nn
&\quad - (-1)^{|u||w|}x_0^{-1}\delta\left(\frac{-x_{2}+x_1}{x_0}\right)
%\cdot\nn
%&\quad\quad\cdot 
(Y^{g})_{WV}^{W}(w, x_{2})
Y_{V}\left(\left({\displaystyle \frac{-x_{2}+x_1}{x_{0}}}\right)^{\mathcal{L}_{g}}
u, x_1\right)v\nn
&= x_1^{-1}\delta\left(\frac{ x_2 + x_0}{x_1}\right)
(Y^{g})_{WV}^{W}\left(Y_{W}^{g}\left(u, x_0\right)w, x_{2}\right)v.
\end{align}
\end{thm}
\pf
This Jacobi identity can be proved by directly using the duality property in Theorem \ref{twisted-io}.
Here we give a proof using  the Jacobi identity (\ref{jacobi-1}).

For $u\in V^{[\alpha]}$ and $v\in V$, 
(\ref{jacobi-.5}) holds. 
Using the definition of $(Y^{g})_{WV}^{W}$ and (\ref{jacobi-.5}), we obtain
\begin{align}\label{jacobi-2}
x_0^{-1}&\delta\left(\frac{x_1 + x_2}{x_0}\right)
\left(\frac{x_{1}+x_{2}}{x_{0}}\right)^{-\alpha}
Y_{W}^{g}\left(\left(
\frac{x_{1}+x_{2}}{x_{0}}\right)^{-\mathcal{N}_{g}}u, x_{0}\right)
(Y^{g})_{WV}^{W}(w, x_{2})v\nn
&\quad - x_0^{-1}\delta\left(\frac{ x_2 + x_1}{x_0}\right)
(Y^{g})_{WV}^{W}(Y_{W}^{g}(u, x_1)w, x_{2})v\nn
&=x_0^{-1}\delta\left(\frac{x_1 + x_2}{x_0}\right)
\left(\frac{x_{1}+x_{2}}{x_{0}}\right)^{-\alpha}
(Y_{W}^{g})_{0}\left(x_{0}^{-\mathcal{N}_{g}}\left(\frac{x_{1}+x_{2}}{x_{0}}\right)^{-\mathcal{N}_{g}}u, x_{0}\right)
(Y^{g})_{WV}^{W}(w, x_{2})v\nn
&\quad - x_0^{-1}\delta\left(\frac{ x_2 + x_1}{x_0}\right)
(Y^{g})_{WV}^{W}(Y_{W}^{g}(u, x_1)w, x_{2})v\nn
&=x_0^{-1}\delta\left(\frac{x_1 + x_2}{x_0}\right)
(Y_{W}^{g})_{0}((x_{1}+x_{2})^{-\mathcal{N}_{g}}u,x_{1}+x_{2})
(Y^{g})_{WV}^{W}(w, x_{2})v\nn
&\quad - x_0^{-1}\delta\left(\frac{ x_2 + x_1}{x_0}\right)
(Y^{g})_{WV}^{W}(Y_{W}^{g}(u, x_1)w, x_{2})v\nn
&=x_0^{-1}\delta\left(\frac{x_1 + x_2}{x_0}\right)
Y_{W}^{g}(u, x_{1}+x_{2})
(Y^{g})_{WV}^{W}(w, x_{2})v\nn
&\quad - x_0^{-1}\delta\left(\frac{ x_2 + x_1}{x_0}\right)
(Y^{g})_{WV}^{W}(Y_{W}^{g}(u, x_1)w, x_{2})v\nn
&=(-1)^{|v||w|}x_0^{-1}\delta\left(\frac{x_1 + x_2}{x_0}\right)
Y_{W}^{g}(u, x_1+x_{2})e^{x_{2}L_{W}^{g}(-1)}
Y_{W}^{g}(v, y)w\lbar_{y^{n}=e^{\pi ni}x^{n},\; \log y=\log x+\pi i}\nn
&\quad - (-1)^{(|u|+|w|)|v|}x_0^{-1}\delta\left(\frac{ x_2 + x_1}{x_0}\right)
e^{x_{2}L_{W}^{g}(-1)}Y_{W}^{g}(v, y)
Y_{W}^{g}(u, x_1)w\lbar_{y^{n}=e^{\pi ni}x^{n},\; \log y=\log x+\pi i}\nn
&=(-1)^{|v||w|}e^{x_{2}L_{W}^{g}(-1)}\left(x_0^{-1}\delta\left(\frac{x_1 -y}{x_0}\right)
Y_{W}^{g}(u, x_1)
Y_{W}^{g}(v, y)w\right.\nn
&\quad\quad\quad\quad \quad\quad\left. - (-1)^{|u||v|}x_0^{-1}\delta\left(\frac{-y + x_1}{x_0}\right)
Y_{W}^{g}(v, y)
Y_{W}^{g}(u, x_1)w\right)\lbar_{y^{n}=e^{\pi ni}x^{n},\; \log y=\log x+\pi i}\nn
&=(-1)^{|v||w|}x_1^{-1}\delta\left(\frac{y+x_0}{x_1}\right)
\left(\frac{y+x_0}{x_1}\right)^{\alpha}e^{x_{2}L_{W}^{g}(-1)}\cdot\nn
&\quad\quad\quad\quad\quad\quad\cdot Y^{g}_{W}\left(Y_{V}\left(\left(\frac{y+x_0}{x_{1}}\right)^{\mathcal{N}_{g}}u, x_0\right)v, y\right)w
\lbar_{y^{n}=e^{\pi ni}x^{n},\; \log y=\log x+\pi i}\nn
&=(-1)^{|u||w|} x_1^{-1}\delta\left(\frac{-x_{2}+x_0}{x_1}\right)
\left(\frac{-x_{2}+x_0}{x_1}\right)^{\alpha} 
\cdot\nn
&\quad\quad\quad\quad\quad\quad\cdot 
(Y^{g})_{WV}^{W}(w, x_{2})
Y_{V}\left(
\left(\frac{-x_{2}+x_0}{x_{1}}\right)^{\mathcal{N}_{g}}
u, x_0\right)v,
\end{align}
where we have used our convention (\ref{convention}).

From (\ref{jacobi-2}), we obtain 
\begin{align}\label{jacobi-3}
x_0^{-1}&\delta\left(\frac{x_1- x_2}{x_0}\right)
\left(\frac{x_{0}+x_{2}}{x_{1}}\right)^{-\alpha}
Y_{W}^{g}\left(\left(\frac{x_{0}+x_{2}}{x_{1}}\right)
^{-\mathcal{N}_{g}}u, x_{1}\right)
(Y^{g})_{WV}^{W}(w, x_{2})v\nn
&\quad - (-1)^{|u||w|} x_0^{-1}\delta\left(\frac{-x_{2}+x_1}{x_0}\right)
\left(\frac{-x_{2}+x_1}{x_0}\right)^{\alpha}  \cdot\nn
&\quad\quad\cdot (Y^{g})_{WV}^{W}(w, x_{2})
Y_{V}\left(
\left(\frac{-x_{2}+x_1}{x_{0}}\right)^{\mathcal{N}_{g}}
u, x_1\right)v\nn
&= x_1^{-1}\delta\left(\frac{ x_2 + x_0}{x_1}\right)
(Y^{g})_{WV}^{W}\left(Y_{W}^{g}\left(u, x_0\right)w, x_{2}\right)v.
\end{align}
Using the property of the formal $\delta$-function, we see that 
the first term in the left-hand side of (\ref{jacobi-3}) is equal to
\begin{align}\label{jacobi-4}
x_0^{-1}&\delta\left(\frac{x_1- x_2}{x_0}\right)
x_{0}^{-\alpha}\left(\frac{\displaystyle 1+\frac{x_{2}}{x_{1}-x_{2}}}{x_{1}}\right)^{-\alpha}\cdot\nn
&\quad\cdot Y_{W}^{g}\left(x_{0}^{-\mathcal{N}_{g}}\left(\frac{\displaystyle 1+\frac{x_{2}}{x_{1}-x_{2}}}{x_{1}}\right)
^{-\mathcal{N}_{g}}u, x_{1}\right)
(Y^{g})_{WV}^{W}(w, x_{2})v.
\end{align}

For $a\in \C$, we have 
\begin{align}\label{formal-idty-1}
\left(\frac{\displaystyle 1+\frac{x_{2}}{x_{1}-x_{2}}}{x_{1}}\right)^{-\alpha}
&= (x_{1}-x_{2})^{\alpha}(x_{1}-x_{2})^{-\alpha}
\left(\frac{\displaystyle 1+\frac{x_{2}}{x_{1}-x_{2}}}{x_{1}}\right)^{-\alpha}\nn
&=(x_{1}-x_{2})^{\alpha}\left(1-\frac{x_{2}}{x_{1}}\right)^{-\alpha}
\left(1+\frac{x_{2}}{x_{1}-x_{2}}\right)^{-\alpha}\nn
&=(x_{1}-x_{2})^{\alpha}\left(\left(1-\frac{x_{2}}{x_{1}}\right)
\left(1+\frac{x_{2}}{x_{1}}\left(1-\frac{x_{2}}{x_{1}}\right)^{-1}\right)\right)^{-\alpha}\nn
& =(x_{1}-x_{2})^{\alpha}.
\end{align}
Similarly, we have
\begin{align*}
\left(\frac{\displaystyle 1+\frac{x_{2}}{x_{1}-x_{2}}}{x_{1}}\right)^{-\mathcal{N}_{g}}
&=x^{\mathcal{N}_{g}}\left(1+\frac{x_{2}}{x_{1}-x_{2}}\right)^{-\mathcal{N}_{g}}\nn
&=x^{\mathcal{N}_{g}}\left(1-\frac{x_{2}}{x_{1}}\right)^{\mathcal{N}_{g}}
\left(1-\frac{x_{2}}{x_{1}}\right)^{-\mathcal{N}_{g}}\left(1+\frac{x_{2}}{x_{1}}\left(1-\frac{x_{2}}{x_{1}}\right)^{-1}\right)^{-\mathcal{N}_{g}}\nn
\end{align*}
\begin{align}\label{formal-idty-2}
&\quad\quad\quad\quad\quad\quad\;\;
=(x_{1}-x_{2})^{\mathcal{N}_{g}}\left(\left(1-\frac{x_{2}}{x_{1}}\right)
\left(1+\frac{x_{2}}{x_{1}}\left(1-\frac{x_{2}}{x_{1}}\right)^{-1}\right)\right)^{-\mathcal{N}_{g}}\nn
&\quad\quad\quad\quad\quad\quad\;\;
=(x_{1}-x_{2})^{\mathcal{N}_{g}}.
\end{align}
Using (\ref{jacobi-4}), 
(\ref{formal-idty-1}) and (\ref{formal-idty-2}), we see that (\ref{jacobi-3}) becomes
\begin{align}\label{jacobi-5}
x_0^{-1}&\delta\left(\frac{x_1- x_2}{x_0}\right)
\left(\frac{x_{1}-x_{2}}{x_{0}}\right)^{\alpha}
Y_{W}^{g}\left(\left(\frac{x_{1}-x_2}{x_0}\right)
^{\mathcal{N}_{g}}u, x_{1}\right)
(Y^{g})_{WV}^{W}(w, x_{2})v\nn
&\quad - (-1)^{|u||w|}x_0^{-1}\delta\left(\frac{-x_{2}+x_1}{x_0}\right)
\left(\frac{-x_{2}+x_1}{x_0}\right)^{\alpha} \cdot\nn
&\quad\quad\cdot (Y^{g})_{WV}^{W}(w, x_{2})
Y_{V}\left(
\left(\frac{-x_{2}+x_1}{x_{0}}\right)^{\mathcal{N}_{g}}
u, x_1\right)v\nn
&= x_1^{-1}\delta\left(\frac{ x_2 + x_0}{x_1}\right)
(Y^{g})_{WV}^{W}(Y_{W}^{g}(u, x_0)w, x_{2})v.
\end{align}
Using our convention (\ref{convention}), we see that
the identity (\ref{jacobi-5}) is (\ref{jacobi-twist}) in the case that $u\in V^{[\alpha]}$.
\epfv

There are many important and useful consequences of the Jacobi identity. For example, 
by taking $\res_{x_{1}}$ on both sides of the Jacobi identity (\ref{jacobi-twist}), we obtain 
an iterate formula for $(Y^{g})_{WV}^{W}$. From this iterate formula, 
we  obtain another generalized weak associativity formula involving the iterate 
$(Y^{g})_{WV}^{W}((Y^{g})_{WV}^{W}(w, x_{0})u, x_{2})v$. 
We shall not give all such consequences in this paper. We give only the following 
generalized commutator formula and generalized weak commutativity which will be used as the main assumptions 
on twist fields in our construction of lower-bounded 
generalized twisted modules in \cite{H-const-twisted-mod} and will play an important role 
in our proof below of the convergence of products of more than one twisted vertex operators and vertex operators for $V$
and one twist vertex operator:

\begin{cor}
For $u\in V, v\in V$ and $w\in W$, we have:
\begin{enumerate}

\item The generalized commutator formula: 
\begin{align}\label{gen-comm-formula}
Y_{W}^{g}&((x_{1}-x_2)
^{\mathcal{L}_{g}}u, x_{1})
(Y^{g})_{WV}^{W}(w, x_{2})v\nn
&\quad - (-1)^{|u||w|}
(Y^{g})_{WV}^{W}(w, x_{2})
Y_{V}((-x_{2}+x_1)^{\mathcal{L}_{g}}
u, x_1)v\nn
&= \res_{x_{0}}x_{0}^{\alpha}x_1^{-1}\delta\left(\frac{ x_2 + x_0}{x_1}\right)
(Y^{g})_{WV}^{W}((Y_{W}^{g})_{0}(u, x_0)w, x_{2})v.
\end{align}

\item The generalized weak commutativity: For $M_{u, w}\in \Z_{+}$  such that  
$x^{\alpha+M_{u, w}}(Y_{W}^{g})_{0}(u, x)w\in W[[x]]$, 
\begin{align}\label{gen-weak-comm}
(x_{1}&-x_{2})^{M_{u, w}}Y_{W}^{g}((x_{1}-x_2)
^{\mathcal{S}_{g}+\mathcal{N}_{g}}u, x_{1})
(Y^{g})_{WV}^{W}(w, x_{2})v\nn
&= (-1)^{|u||w|} (x_{1}-x_{2})^{M_{u, w}}
(Y^{g})_{WV}^{W}(w, x_{2})
Y_{V}((-x_{2}+x_1)^{\mathcal{S}_{g}+\mathcal{N}_{g}}
u, x_1)v.
\end{align}
\end{enumerate}
\end{cor}
\pf
Let $u\in V^{[\alpha]}$. Then the Jacobi identity (\ref{jacobi-twist}) becomes 
(\ref{jacobi-5}). Multiplying both sides of (\ref{jacobi-5}) by $x_{0}^{\alpha+m}$ for $m\in \Z_{+}$, 
we see that the two sides contains only integral powers of $x_{0}$. Then applying $\res_{x_{0}}$
to both sides and using (\ref{y-g-w-0}), we obtain
\begin{align}\label{gen-comm-formula-1}
(x_{1}&-x_{2})^{\alpha+m}
Y_{W}^{g}\left(\left(\frac{x_{1}-x_2}{x_0}\right)
^{\mathcal{N}_{g}}u, x_{1}\right)
(Y^{g})_{WV}^{W}(w, x_{2})v\nn
&\quad - (-1)^{|u||w|}(-x_{2}+x_1)^{\alpha+m} (Y^{g})_{WV}^{W}(w, x_{2})
Y_{V}\left(
\left(\frac{-x_{2}+x_1}{x_{0}}\right)^{\mathcal{N}_{g}}
u, x_1\right)v\nn
&= \res_{x_{0}}x_{0}^{\alpha+m}x_1^{-1}\delta\left(\frac{ x_2 + x_0}{x_1}\right)
(Y^{g})_{WV}^{W}(Y_{W}^{g}(u, x_0)w, x_{2})v\nn
&= \res_{x_{0}}x_{0}^{\alpha+m}x_1^{-1}\delta\left(\frac{ x_2 + x_0}{x_1}\right)
(Y^{g})_{WV}^{W}((Y_{W}^{g})_{0}(x_{0}^{-\mathcal{N}_{g}}u, x_0)w, x_{2})v.
\end{align}
Taking $u$ in (\ref{gen-comm-formula-1}) to be
 $x_{0}^{\mathcal{N}_{g}}u$, we
see that (\ref{gen-comm-formula-1}) becomes 
\begin{align}\label{gen-comm-formula-2}
(x_{1}&-x_{2})^{\alpha+m}
Y_{W}^{g}((x_{1}-x_2)
^{\mathcal{N}_{g}}u, x_{1})
(Y^{g})_{WV}^{W}(w, x_{2})v\nn
&\quad - (-1)^{|u||w|}(-x_{2}+x_1)^{\alpha+m} (Y^{g})_{WV}^{W}(w, x_{2})
Y_{V}((-x_{2}+x_1)^{\mathcal{N}_{g}}
u, x_1)v\nn
&= \res_{x_{0}}x_{0}^{\alpha+m}x_1^{-1}\delta\left(\frac{ x_2 + x_0}{x_1}\right)
(Y^{g})_{WV}^{W}((Y_{W}^{g})_{0}(u, x_0)w, x_{2})v.
\end{align}

Taking $m=0$ in (\ref{gen-comm-formula-2}), we obtain 
(\ref{gen-comm-formula}) in the case $u\in V^{[\alpha]}$.

Since
$x_{0}^{\alpha+M_{u, w}}(Y_{W}^{g})_{0}(u, x_0)w$
is a power series in $x_{0}$, for $m=M_{u, w}$, the right-hand side of 
(\ref{gen-comm-formula-2}) is equal to $0$. Thus the left-hand side of
(\ref{gen-comm-formula-2}) is equal to $0$ when $m=M_{u, w}$, proving 
(\ref{gen-weak-comm}). 
\epfv

The generalized commutator formula and 
generalized weak commutativity above are formulated in terms of $Y_{W}^{g}((x_{1}-x_2)
^{\mathcal{L}_{g}}u, x_{1})$ and 
$Y_{V}((-x_{2}+x_1)^{\mathcal{L}_{g}}
u, x_1)$. 
In the construction of lower-bounded generalized 
twisted modules in \cite{H-const-twisted-mod}, we 
need a generalized commutator formula 
and generalized weak commutativity
expressed in terms of twisted fields and untwisted fields without $(x_{1}-x_2)
^{\mathcal{L}_{g}}$ and $(-x_{2}+x_1)^{\mathcal{L}_{g}}$. 
Here we rewrite the generalized commutator formula and 
generalized weak commutativity above as follows:

\begin{cor}
For $u\in V^{[\alpha]}, v\in V$ and $w\in W$, the generalized 
commutator formula (\ref{gen-comm-formula}) and generalized weak commutativity 
(\ref{gen-comm-formula}) can be rewritten as 
\begin{align}\label{gen-comm-formula-3}
(x&_{1}-x_{2})^{\alpha}(x_{1}-x_2)
^{\mathcal{N}_{g}}Y_{W}^{g}(u, x_{1})(x_{1}-x_2)
^{-\mathcal{N}_{g}}
(Y^{g})_{WV}^{W}(w, x_{2})v\nn
&\quad - (-1)^{|u||w|} (-x_{2}+x_{1})^{\alpha}
(Y^{g})_{WV}^{W}(w, x_{2})(-x_{2}+x_1)^{\mathcal{N}_{g}}
Y_{V}(u, x_1)
(-x_{2}+x_1)^{-\mathcal{N}_{g}}v\nn
&=  \sum_{k=0}^{M_{u, w}-1}\frac{1}{k!}
x_1^{-1}\frac{\partial^{k}}{\partial x_{2}^{k}}\delta\left(\frac{ x_2}{x_1}\right)
(Y^{g})_{WV}^{W}((Y_{W}^{g})_{\alpha+k, 0}(u)w, x_{2})v.
\end{align}
and
\begin{align}\label{gen-weak-comm-3}
(x_{1}&-x_{2})^{\alpha+M_{u, w}} (x_{1}-x_{2})^{\mathcal{N}_{g}}Y_{W}^{g}(u, x_{1}) 
(x_{1}-x_{2})^{-\mathcal{N}_{g}}
(Y^{g})_{WV}^{W}(w, x_{2})v\nn
&\;\;=(-1)^{|u||w|} (-x_{2}+x_{1})^{\alpha+M_{u, w}}
(Y^{g})_{WV}^{W}(w, x_{2})\cdot\nn
&\quad\quad\quad\quad\quad\quad\cdot 
(-x_{2}+x_1)^{\mathcal{N}_{g}} Y_{V}(u, x_1)
(-x_{2}+x_1)^{-\mathcal{N}_{g}}v,
\end{align} 
respectively. In particular, when $g$ is semisimple (for example, when $g$ is of finite order), 
we have 
\begin{align}\label{gen-comm-formula-4}
(x&_{1}-x_{2})^{\alpha}Y_{W}^{g}(u, x_{1})
(Y^{g})_{WV}^{W}(w, x_{2})v\nn
&\quad - (-1)^{|u||w|} (-x_{2}+x_{1})^{\alpha}
(Y^{g})_{WV}^{W}(w, x_{2})
Y_{V}(u, x_1)v\nn
&= \sum_{k=0}^{M_{u, w}-1}\frac{1}{k!}
x_1^{-1}\frac{\partial^{k}}{\partial x_{2}^{k}}\delta\left(\frac{ x_2}{x_1}\right)
(Y^{g})_{WV}^{W}((Y_{W}^{g})_{\alpha+k, 0}(u)w, x_{2})v.
\end{align}
and
\begin{align}\label{gen-weak-comm-4}
(x_{1}&-x_{2})^{\alpha+M_{u, w}} Y_{W}^{g}(u, x_{1}) 
(Y^{g})_{WV}^{W}(w, x_{2})v\nn
&\;\;=(-1)^{|u||w|} (-x_{2}+x_{1})^{\alpha+M_{u, w}}
(Y^{g})_{WV}^{W}(w, x_{2})
Y_{V}(u, x_1)v,
\end{align} 
\end{cor}
\pf
In the case $u\in V^{[\alpha]}$, 
\begin{align*}
(x_{1}-x_{2})^{\mathcal{S}_{g}}u&=(x_{1}-x_{2})^{\alpha},\nn
(-x_{2}+x_1)^{\mathcal{S}_{g}}u&=
(-x_{2}+x_{1})^{\alpha}.
\end{align*}
Also
\begin{align*}
(x_{1}-x_2)^{\mathcal{N}_{g}}u&=\sum_{k\in \N}\frac{1}{k!}
(\log (x_{1}-x_{2}))^{k}\mathcal{N}_{g}^{k}u,\nn
(-x_{2}+x_1)^{\mathcal{N}_{g}}u&=\sum_{k\in \N}\frac{1}{k!}
(\log (-x_{2}+x_1))^{k}\mathcal{N}_{g}^{k}u.
\end{align*}
Using these formulas, (\ref{V-formal-N-g-conj}), (\ref{formal-N-g-conj}),
$(Y_{W}^{g})_{0}(u, x_0)w=\sum_{n\in \alpha+\Z}(Y_{W}^{g})_{n, 0}(u)wx_{0}^{-n-1}$
and the definition of $M_{u, w}$, 
we see that (\ref{gen-comm-formula}) and 
(\ref{gen-weak-comm}) become (\ref{gen-comm-formula-3}) and
(\ref{gen-weak-comm-3}), respectively.
\epfv

We now study products of more than one twisted vertex operators or vertex operators for $V$
and one 
twist
vertex operators:

\begin{thm}\label{k+l-prod-twted-twt-vo}
For $w'\in W'$, $v_{1}\in V^{[\alpha_{1}]}, \dots, v_{k+l}\in V^{[\alpha_{k+l}]}, 
v\in V^{[\alpha]}$ and $w\in W$, 
the series
\begin{equation}\label{k+l-prod-twted-twt-vo-1}
\langle w',  (Y_{W}^{g})^{p}(v_{1}, z_{1})\cdots (Y_{W}^{g})^{p}(v_{k}, z_{k})
((Y^{g})_{WV}^{W})^{p}(w, z)
Y_{V}(v_{k+1}, z_{k+1})\cdots Y_{V}(v_{k+l}, z_{k+l})v\rangle
\end{equation}
is absolutely convergent in the region $|z_{1}|>\cdots>|z_{k}|>|z|>
|z_{k+1}|>\cdots >|z_{k+l}|>0$. Moreover, there exists a 
multivalued analytic function of the form 
\begin{align}\label{k+l-prod-twted-twt-correl-1}
&\sum_{n_{1}, \dots, n_{k+l}, n=0}^{N}f_{n_{1}\cdots n_{k+l}n}(z_{1}, \dots,
z_{k+l}, z)\cdot\nn
&\quad\quad\quad\quad \cdot  
(z_{1}-z)^{-\alpha_{1}}\cdots (z_{k+l}-z)^{-\alpha_{k+l}}z^{-\alpha}
(\log (z_{1}-z))^{n_{1}}\cdots (\log (z_{k+l}-z))^{n_{k+l}}(\log z)^{n},
\end{align}
denoted by 
$$F(\langle w',  Y_{W}^{g}(v_{1}, z_{1})\cdots Y_{W}^{g}(v_{k}, z_{k})
(Y^{g})_{WV}^{W}(w, z)
Y_{V}(v_{k+1}, z_{k+1})\cdots Y_{V}(v_{k+l}, z_{k+l})v\rangle),$$
where $N\in \N$ and $f_{n_{1}\cdots n_{k+l}n}(z_{1}, \dots,
z_{k+l}, z)$ for $n_{1}, \cdots, n_{k+l}, n=0, \dots, N$ are rational functions
of $z_{1}, \dots, z_{k+l}, z$ with the only possible poles $z_{i}=0$ for
$i=1, \dots, k+l$, $z=0$, $z_{i}-z_{j}=0$
for $i, j=1, \dots, k+l$, $i\ne j$, $z_{i}-z=0$ for $i=1, \dots, k+l$, 
such that the sum of (\ref{k+l-prod-twted-twt-vo-1})  is equal to the branch
\begin{align}\label{k+l-prod-twted-twt-correl-2}
F&^{p}(\langle w',  Y_{W}^{g}(v_{1}, z_{1})\cdots Y_{W}^{g}(v_{k}, z_{k})
(Y^{g})_{WV}^{W}(w, z)
Y_{V}(v_{k+1}, z_{k+1})\cdots Y_{V}(v_{k+l}, z_{k+l})v\rangle)\nn
&=\sum_{n_{1}, \dots, n_{k+l}, n=0}^{N}f_{n_{1}\cdots n_{k+l}n}(z_{1}, \dots,
z_{k+l}, z)\cdot\nn
&\quad\quad\quad\quad\quad \cdot  
e^{-\alpha_{1}l_{p}(z_{1}-z)}\cdots e^{-\alpha_{k+l}l_{p}(z_{k+l}-z)}
e^{-\alpha l_{p}(z)}
(l_{p}(z_{1}-z))^{n_{1}}\cdots (l_{p}(z_{k+l}-z))^{n_{k+l}}(l_{p}(z))^{n},
\end{align}
of (\ref{k+l-prod-twted-twt-correl-1})
in the region 
given by $|z_{1}|>\cdots>|z_{k}|>|z|>
|z_{k+1}|>\cdots >|z_{k+l}|>0$, $|\arg (z_{i}-z)-\arg z_{i}|<\frac{\pi}{2}$ for 
$i=1, \dots, k$ and $|\arg (z_{i}-z)-\arg z|<\frac{\pi}{2}$ for 
$i=k+1, \dots, k+l$. 
In addition,  the orders of the pole $z_{j}=0$ of the rational functions 
$f_{n_{1}\cdots n_{k+l}n}(z_{1}, \dots,z_{k+l}, z)$ 
have a lower bound independent of $v_{q}$ for $q\ne j$, $w$ and $w'$; 
the orders of the pole $z=0$ of the rational functions 
$f_{n_{1}\cdots n_{k+l}n}(z_{1}, \dots, z_{k+l}, z)$ 
have a lower bound independent of $v_{1}, \dots, v_{k+l}$ and $w'$; 
the orders of the pole $z_{j}=z_{m}$ of the rational functions 
$f_{n_{1}\cdots n_{k+l}n}(z_{1}, \dots, z_{k+l}, z)$ 
have a lower bound independent of $v_{q}$ for $q\ne j, m$, $v$, $w$ and $w'$;
the orders of the pole $z_{j}=z$ of the rational functions 
$f_{n_{1}\cdots n_{k+l}n}(z_{1}, \dots, z_{k+l}, z)$ 
have a lower bound independent of $v_{q}$ for $q\ne j$, $v$ and $w'$.
\end{thm}
\pf
Let $L_{W}(-1)'$ be the adjoint of $L_{W}(-1)$ on $W'$. By Theorem \ref{k-prod-twted-vo}, 
\begin{align}
\langle & e^{zL_{W}(-1)'}w',  (Y_{W}^{g})^{p}(v_{1}, z_{1}-z)
\cdots (Y_{W}^{g})^{p}(v_{k}, z_{k}-z)\cdot\nn
&\quad\quad\quad\quad \quad\quad\cdot 
(Y_{W}^{g})^{p}(v_{k+1}, z_{k+1}-z)\cdots (Y_{W}^{g})^{p}(v_{k+l}, z_{k+l}-z)
(Y_{W}^{g})^{p}(v, -z)w\rangle
\end{align}
converges absolutely in the region $|z_{1}-z|>\cdots >|z_{k+l}-z|>|z|>0$
to 
\begin{align}\label{k+l-prod-twted-twt-correl-3}
&\sum_{n_{1}, \dots, n_{k}, n=0}^{N}g_{n_{1}\cdots n_{k+l}n}(z_{1}-z, \dots,
z_{k+l}-z, -z)\cdot\nn
&\quad\quad\quad\quad \cdot 
e^{-\alpha_{1}l_{p}(z_{1}-z)}\cdots e^{-\alpha_{k+l}l_{p}(z_{k+l}-z)}
e^{-\alpha l_{p}(-z)}
(l_{p}(z_{1}-z))^{n_{1}}\cdots (l_{p}(z_{k+l}-z))^{n_{k+l}}(l_{p}(z))^{n},
\end{align}
where $n_{1}, \dots, n_{k+l}\in \N$
and $g_{n_{1}\cdots n_{k+l}n}(z_{1}-z, \dots,
z_{k+l}-z, -z)$ are rational functions in $z_{1}-z, \dots,
z_{k+l}-z, -z$ with the only possible poles at $z_{i}-z=0$, $z=0$ and $z_{i}-z_{j}=0$
for $i\ne j$. Moreover,  the orders of the pole $z_{j}=0$ of the rational functions 
$g_{n_{1}\cdots n_{k+l}n}(z_{1}-z, \dots,
z_{k+l}-z, -z)$
have a lower bound independent of $v_{q}$ for $q\ne j$, $w$ and $w'$; 
the orders of the pole $z=0$ of the rational functions 
$g_{n_{1}\cdots n_{k+l}n}(z_{1}-z, \dots,
z_{k+l}-z, -z)$ 
have a lower bound independent of $v_{1}, \dots, v_{k+l}$ and $w'$; 
the orders of the pole $z_{j}=z_{m}$ of the rational functions 
$g_{n_{1}\cdots n_{k+l}n}(z_{1}-z, \dots,
z_{k+l}-z, -z)$
have a lower bound independent of $v_{q}$ for $q\ne j, m$, $v$, $w$ and $w'$;
the orders of the pole $z_{j}=z$ of the rational functions 
$g_{n_{1}\cdots n_{k+l}n}(z_{1}-z, \dots,
z_{k+l}-z, -z)$
have a lower bound independent of $v_{q}$ for $q\ne j$, $v$ and $w'$.
Let 
$$f_{n_{1}\cdots n_{k+l}n}(z_{1}, \dots,
z_{k+l}, z)=g_{n_{1}\cdots n_{k+l}n}(z_{1}-z, \dots,
z_{k+l}-z, -z)e^{-\alpha\pi i}.$$
Since $l_{p}(-z)=l_{p}(z)+\pi i$ when $0\le \arg z<\pi$, 
(\ref{k+l-prod-twted-twt-correl-3}) is equal to a multivalued analytic function 
of the form 
(\ref{k+l-prod-twted-twt-correl-1}) with branches of the form 
(\ref{k+l-prod-twted-twt-correl-2}) satisfying all the properties
when $0\le \arg z<\pi$. Using our notations, we denote this function 
by 
\begin{align*}
F(\langle & e^{zL_{W}(-1)'}w',  Y_{W}^{g}(v_{1}, z_{1}-z)
\cdots Y_{W}^{g}(v_{k}, z_{k}-z)\cdot\nn
&\quad\quad\quad\quad \quad\cdot 
Y_{W}^{g}(v_{k+1}, z_{k+1}-z)\cdots Y_{W}^{g}(v_{k+l}, z_{k+l}-z)
Y_{W}^{g}(v, -z)w\rangle)
\end{align*}

On the other hand, when $0\le \arg z<\pi$, 
using the $L(-1)$-commutator formula, the $L(-1)$-derivative 
property and the associativity for the twisted vertex operators, we have
\begin{align}\label{k+l-prod-twted-twt-correl-4}
F^{p}(\langle & e^{zL_{W}(-1)'}w',  Y_{W}^{g}(v_{1}, z_{1}-z)
\cdots Y_{W}^{g}(v_{k}, z_{k}-z)\cdot\nn
&\quad\quad\quad\quad \quad\cdot 
Y_{W}^{g}(v_{k+1}, z_{k+1}-z)\cdots Y_{W}^{g}(v_{k+l}, z_{k+l}-z)
Y_{W}^{g}(v, -z)w\rangle)\nn
&=F^{p}(\langle w',  Y_{W}^{g}(v_{1}, z_{1})
\cdots Y_{W}^{g}(v_{k}, z_{k})\cdot\nn
&\quad\quad\quad\quad \quad\cdot e^{zL_{W}(-1)}
Y_{W}^{g}(v_{k+1}, z_{k+1}-z)\cdots Y_{W}^{g}(v_{k+l}, z_{k+l}-z)
Y_{W}^{g}(v, -z)w\rangle)\nn
&=F^{p}(\langle w',  Y_{W}^{g}(v_{1}, z_{1})
\cdots Y_{W}^{g}(v_{k}, z_{k})\cdot\nn
&\quad\quad\quad\quad \quad\cdot e^{zL_{W}(-1)}
Y_{W}^{g}(Y_{V}(v_{k+1}, z_{k+1})\cdots Y_{V}(v_{k+l}, z_{k+l})v, -z)w\rangle)\nn
&=F^{p}(\langle  w',   Y_{W}^{g}(v_{1}, z_{1})\cdots Y_{W}^{g}(v_{k}, z_{k})
(Y^{g})_{WV}^{W}(w, z)
Y_{V}(v_{k+1}, z_{k+1})\cdots Y_{V}(v_{k+l}, z_{k+l})v\rangle).
\end{align}
Then by analytic extensions, (\ref{k+l-prod-twted-twt-correl-4}) holds without the condition 
$0\le \arg z<\pi$. 
But (\ref{k+l-prod-twted-twt-vo-1}) is a series of the same form as the expansion 
of the left-hand side of (\ref{k+l-prod-twted-twt-correl-4}) in the region 
given by $|z_{1}|>\cdots>|z_{k}|>|z|>
|z_{k+1}|>\cdots >|z_{k+l}|>0$, $|\arg (z_{i}-z)-\arg z_{i}|<\frac{\pi}{2}$ for 
$i=1, \dots, k$ and $|\arg (z_{i}-z)-\arg z|<\frac{\pi}{2}$ for 
$i=k+1, \dots, k+l$.  By (\ref{k+l-prod-twted-twt-correl-4}),
(\ref{k+l-prod-twted-twt-vo-1}) must be the expansion of the
left-hand side of (\ref{k+l-prod-twted-twt-correl-4}).
Since  we have proved that the left-hand side of (\ref{k+l-prod-twted-twt-correl-4})
is of the form (\ref{k+l-prod-twted-twt-correl-2}) satisfying all the properties,
the result is proved.
\epf

\begin{cor}\label{k-comm}
For $v_{1}, \dots,  v_{k-1}\in V$, $v\in V$, 
$w\in W$, $w'\in W'$,
$p\in \Z$, $\tau\in S_{k}$ and fixed $1\le i\le k$, 
$$
F^{p}(\langle w', \varphi_{1}(z_{1})\cdots \varphi_{k}(z_{k})v\rangle)
=\pm F^{p}(\langle w', \varphi_{\tau(1)}(z_{\tau(1)})\cdots
\varphi_{\tau(k)}(z_{\tau(k)})v\rangle),
$$
where $\varphi_{j}(z_{j})=Y_{W}^{g}(v_{j}, z_{j})$
for $j\ne i$ and $\varphi_{i}=(Y^{g})_{WV}^{W}(w, z_{l})$ and
the sign $\pm$ is uniquely determined by $\tau$
and $|v_{1}|, \dots, |v_{k-1}|, |w|$.
\end{cor}
\pf
This result follows immediately from 
Theorem \ref{k-prod-twted-vo}, the duality property for $Y_{W}^{g}$
and Corollary \ref{twisted-twist-comm}.
\epfv

\noindent {\small \sc Department of Mathematics, Rutgers University,
110 Frelinghuysen Rd., Piscataway, NJ 08854-8019}

\noindent {\em E-mail address}: yzhuang@math.rutgers.edu

\end{document}